\documentclass[11pt]{article}
\usepackage[T1]{fontenc}
\usepackage{lmodern}
\usepackage{amsmath,amssymb,amsthm}
\usepackage[margin=1.2in]{geometry}
\usepackage{microtype}
\usepackage{needspace}
\theoremstyle{plain}
\newtheorem{theorem}{Theorem}[section]
\newtheorem{lemma}[theorem]{Lemma}
\newtheorem{proposition}[theorem]{Proposition}
\newtheorem{corollary}[theorem]{Corollary}
\theoremstyle{definition}
\newtheorem{definition}[theorem]{Definition}
\newtheorem{example}[theorem]{Example}

\theoremstyle{remark}
\newtheorem{remark}[theorem]{Remark}
\newcommand{\Z}{\mathbb{Z}}
\newcommand{\N}{\mathbb{N}}
\newcommand{\R}{\mathbb{R}}
\newcommand{\lam}{\lambda^{-}}
\newcommand{\lams}{\lambda^{-}_{s}}
\newcommand{\D}{D_0}
\pdftrailerid{}
\usepackage[hidelinks]{hyperref}
\hypersetup{
  pdftitle={The Structure of a Non-Isometric Non-Contractive Bijection of a Hereditarily Precompact Quasi-Metric Space},
  pdfauthor={Thomas Duzer}
}

\hypersetup{pdfkeywords={quasi-metric space, expand-contract plasticity,
    non-contractive map, hereditarily precompact, toll spiral,
    monothetic group, paratopological group,
    formal verification, Lean}}

\title{The Structure of a Non-Isometric Non-Contractive\\
Bijection of a Hereditarily Precompact Quasi-Metric Space}
\author{Thomas Duzer\thanks{thomas.duzer@gmail.com.\newline
  ORCID: \href{https://orcid.org/0009-0007-2147-081X}{\mbox{0009-0007-2147-081X}}.}}
\date{}
\begin{document}
\maketitle

\begin{abstract}
Let $F$ be a non-contractive bijection of a hereditarily precompact
quasi-metric space $(X, d)$. We show that $F$ is asymptotically an
isometry along its backward orbits: for every $a \in X$ and
$\varepsilon > 0$, $0 \le d(Fx, Fy) - d(x, y) < \varepsilon$ for
all $y$ and all $x$ sufficiently far in the backward orbit of $a$.
More precisely, $d = \lam + \D$, where the \emph{asymptotic distance}
$\lam(x, y) = \lim_k d(F^{-k}x, F^{-k}y)$ is the largest
quasi-pseudo-metric below $d$ for which $F$ is an isometry, and the
defect $\D(x, y)$ lies between explicit bounds given by $\lam$
and by a function $c$ that is nondecreasing along orbits and tends
to $0$ along backward orbits. On each
orbit, $\lam$ induces, after
completion, an invariant quasi-pseudo-metric on a compact metrizable monothetic
group on which $F$ acts as a rotation, with the compact group
topology.
Conversely, every infinite compact metrizable monothetic group, with
any rotation by a generator and any invariant quasi-pseudo-metric whose
symmetrization is a compatible metric, arises in this way from a map
that is not an isometry. Every such $c$ occurs over every possible
$\lam$, and the defects on a single orbit are characterized. These results give a structural answer to a classification
question from the author's paper on the toll
spiral (arXiv:2610.13973), a non-isometric
bijection of this kind answering a question of
Zavarzina. Three examples on the integers show that the defect is
not determined by $c$, and that the question, read as asking for a
toll charged on descents only, has a negative answer, even up to
bi-Lipschitz conjugacy. The results hold for quasi-pseudo-metric
spaces in which every sequence has a left K-Cauchy
subsequence. Except for the statements listed at
the end, among them the completion to a compact group, the results
are machine-checked in Lean.
\end{abstract}


\needspace{6\baselineskip}
\section{Introduction}\label{sec:intro}

A self-map $F$ of a quasi-metric space $(X, d)$ is \emph{non-contractive}
if $d(Fx, Fy) \ge d(x, y)$ for all $x, y$. On a totally bounded metric
space every non-contractive self-map is an isometry
\cite{FH1936, NPW2006}; this rigidity is the origin of the
expand--contract plasticity literature (see \cite{BMZ2026} and the
references there). Zavarzina \cite{Zav2023} showed that this rigidity fails on
hereditarily precompact quasi-metric spaces, and asked whether it
holds at least for bijections of precompact quasi-metric spaces. The \emph{toll spiral}
\cite{Duzer2026} shows that it does not. It is restored, however,
when the topology of $d$ or that of its conjugate is compact, even
countably compact \cite{DuzerC}. Agyingi
\cite{Agyingi2026} studies nonexpansive bijections (the inverses of
ours) of $T_0$ quasi-pseudo-metric spaces and their $q$-hyperconvex hulls,
and \cite{DuzerC} compares that setting with ours.

\paragraph{Terminology.}
We use the language of the toll-spiral paper \cite{Duzer2026}, in
which $d(x,y)$ is the cost of a trip from $x$ to $y$: a trip back may
cost a different amount, and the triangle inequality says that a
detour is never cheaper than the direct route. The \emph{asymptotic distance}
$\lam(x, y)$ is the limit of the cost of the trip from $x$ to $y$
when both points are moved back along their orbits. The map
preserves it, and the \emph{defect} is the difference $d - \lam$.
The \emph{descent defect} $c(x)$ is the least upper bound of the
defects of the trips from $x$ to earlier points of its orbit.
In the examples
of Section~\ref{sec:spirals}, $d$ is a base metric plus a surcharge,
the \emph{toll}, and the toll is the defect.

Between the counterexample and the compactness theorem lies a
classification question, Question~7.6 of the toll-spiral paper. The
toll spiral is obtained from an isometry, a rotation of the circle
restricted to one orbit, by adding a fixed surcharge to every trip
from a point of a forward-invariant set to an earlier point of its
orbit. The question asks whether every bijective non-contractive
non-isometry of a hereditarily precompact quasi-metric space arises
in this way, up to a suitable conjugacy, possibly with infinitely
many such sets. We show that every such map has a common structure
(Theorem~\ref{thm:structure}), and that this structure still allows
the surcharge to take many forms (Section~\ref{sec:spirals}).

\subsection*{What is forced}

Fix a quasi-pseudo-metric space $(X, d)$ satisfying hypothesis (i)
of \cite{DuzerC}: every sequence has a left K-Cauchy subsequence.
For quasi-metric spaces, this is hereditary precompactness. Fix also
a non-contractive bijection $F$ of $X$, with inverse $G$. The
\emph{asymptotic distance} is
\[
\lam(x, y) = \lim_{k \to \infty} d(G^k x, G^k y);
\]
the \emph{defect} is $\D = d - \lam$, and the \emph{symmetrized
asymptotic distance} is $\lams(x, y) = \lam(x, y) + \lam(y, x)$. In a pair
$(x, y)$ we call $x$ the \emph{source} and $y$ the \emph{target}; the
pair is \emph{strict} if $\D(x, y) > 0$, a \emph{descent} if
$y = G^\Delta x$ and an \emph{ascent} if $y = F^\Delta x$ for some
$\Delta \ge 1$. The \emph{descent defect} of $x$ is the
supremum of the defects of its backward iterates,
$c(x) = \sup_{\Delta \ge 1} \D(x, G^\Delta x)$. The main theorem
gathers what is forced.

\begin{theorem}[Structure]\label{thm:structure}
Let $(X, d)$ satisfy (i) and let $F$ be a non-contractive bijection
of $X$ with inverse $G$. Then:
\begin{enumerate}
\item[(a)] $\lam$ is a quasi-pseudo-metric with $\lam \le d$, for
which $F$ is an isometry, and whose symmetrization is totally
bounded on every orbit;
\item[(b)] for all $x, y$ and $k \ge 0$,
\[
c(G^k x) - c(G^k y) - \lams(x, y) \;\le\; \D(x, y) \;\le\; c(x),
\]
where the descent defect $c$ is nondecreasing along orbits, $c(x) \le c(Fx)$,
and tends to $0$ along backward orbits, $c(G^k x) \to 0$.
\end{enumerate}
Consequently, for every $\varepsilon > 0$, the set
$U_\varepsilon = \{x : c(x) \ge \varepsilon\}$ satisfies
$F(U_\varepsilon) \subseteq U_\varepsilon$, and every backward orbit
eventually leaves it. Every pair
whose source lies outside $U_\varepsilon$ has defect less than
$\varepsilon$: if $x \notin U_\varepsilon$, then
$d(x, y) < \lam(x, y) + \varepsilon$ for every $y$.
\end{theorem}

The double inequality in (b) bounds the
defect of a pair above by the descent defect of its source, and below by the
difference between the descent defects of source and target, now or at any
earlier time, less their symmetrized asymptotic distance.

The theorem's sharpest consequence is simple to state: $F$ is
asymptotically an isometry along its backward orbits. More precisely
(Corollary~\ref{cor:rigid}), for every $a \in X$ and
$\varepsilon > 0$ there is $k$ such that
\[
0 \;\le\; d(Fx, Fy) - d(x, y) \;<\; \varepsilon
\qquad\text{for } x = G^{k'} a,\ k' > k,\ \text{and every } y.
\]
Moreover, every backward orbit is totally bounded for the symmetrized
distance $d(x, y) + d(y, x)$, and for every $\varepsilon > 0$ it has a tail on which some
power $G^n$, with $n$ as large as we like, moves every point by less
than $\varepsilon$ in both directions (see the two paragraphs after
Corollary~\ref{cor:rigid}). With the forward returns
of \cite[Lemma~3.3]{DuzerC}, these near-returns show that a space satisfying (i) that
carries a non-isometric non-contractive bijection is contained in no
$R_0$ space in which its left K-Cauchy sequences have upper limits, and in
none in which they have lower limits; in particular, in no
sequentially Yoneda complete quasi-metric space and no quasi-metric
space that is countably compact for its topology or its conjugate
topology (Corollary~\ref{cor:host}).

Section~\ref{sec:tollshadow} recasts Theorem~\ref{thm:structure} as
a structural answer to \cite[Question~7.6]{Duzer2026}
(Corollary~\ref{cor:tollanswer}).

\subsection*{What is not forced}

Section~\ref{sec:spirals} constructs three examples on $\Z$,
spirals over the circle. Each is a $T_1$ quasi-metric satisfying
(i), on which the shift is non-contractive but not an isometry, and
on each the defect behaves differently within the bounds of
Theorem~\ref{thm:structure}(b). In the height spiral, the defect
of a descent depends on its target and not only on its source. In
the thin spiral, a source with positive descent defect has descents
with no defect (for a suitable value of its parameter). In the
ascent spiral, ascents carry defect. On a single infinite orbit, a
map has the form of the toll spiral, with possibly infinitely many
forward-invariant sets of descents (descent zones), exactly when its
ascents carry no defect (Proposition~\ref{prop:tolled}). So the
ascent spiral answers \cite[Question~7.6]{Duzer2026} negatively in
this form, even up to bi-Lipschitz conjugacy
(Proposition~\ref{prop:bilip}).

Nor is any rate of decay forced. The descent defect always tends to $0$ along
backward orbits, but in the height spiral it is positive at every
point. For a class of two-valued tolls that
includes the toll spiral, the descent defect is zero far enough
back on each orbit (Theorem~\ref{thm:zone}). In fact nothing beyond
Theorem~\ref{thm:structure} constrains the descent defect: for every
quasi-pseudo-metric $\ell$ that satisfies (i) and for which $F$ is an
isometry, every function that is
nondecreasing along orbits and tends to $0$ along backward orbits is
the descent defect of $F$ for some quasi-pseudo-metric $\ell + \tau$
that satisfies (i), for which $F$ is non-contractive and whose
asymptotic distance is $\ell$, and the toll
$\tau_c$ of Section~\ref{sec:band} is the least such defect that
satisfies the triangle inequality (Theorem~\ref{thm:occur}). On a
single orbit, the defects with a given asymptotic distance and
descent defect are characterized (Theorem~\ref{thm:orbit}), and every
pair $(c, A)$ of nondecreasing functions with $0 \le A \le c$ and $c$
tending to $0$ in the past is the pair of descent and ascent defects
of some toll (Theorem~\ref{thm:pairs}). Given
(i) for $\ell$, condition (i) for $\ell + \tau$ is a condition of
Ramsey type: no sequence $(x_i)$ may have $\tau(x_i, x_j) \ge
\varepsilon > 0$ for all $i < j$ (Theorem~\ref{thm:ramsey}).

Nor is the asymptotic distance on an orbit constrained beyond what
Section~\ref{sec:shadow} proves. Take a compact metrizable monothetic
group $K$, an element $g$ generating a dense subgroup, and an invariant
quasi-pseudo-metric on $K$ whose symmetrization is a metric inducing
its topology. Then $K$, with this quasi-pseudo-metric and the rotation
by $g$, is, up to isomorphism, the completion of an orbit for the
asymptotic distance of a non-contractive bijection that is not an
isometry, on a space satisfying (i) (Proposition~\ref{prop:realize}).

\subsection*{Methods and plan}

Section~\ref{sec:prelim} sets up the asymptotic distance, the defect and the
descent defect. The
proof of Theorem~\ref{thm:structure}, assembled at the end of
Section~\ref{sec:depth}, rests on three facts:
\begin{itemize}
\item[(1)] for a strict pair $(a, b)$, along the differences $\Delta$ of a chain, the target
nearly returns to itself, $d(G^\Delta b, b) < \varepsilon$, while
$d(a, G^\Delta a)$ is almost as large as the defect of $(a, b)$; so the
source carries the defect, and the orbit of $a$ is infinite and,
with the restricted distance, itself a non-isometric example
(Section~\ref{sec:source});
\item[(2)] the defect of every pair $(a, y)$ is at least the
difference of the descent defects of $a$ and $y$, now or at any
earlier time, less their symmetrized asymptotic distance, and at most
the descent defect of $a$, which is nondecreasing along orbits
(Section~\ref{sec:band});
\item[(3)] on every orbit the symmetrized asymptotic distance is totally
bounded, so the orbit, with $\lam$, maps densely into a
compact monothetic group on which $F$ acts as a rotation
(Section~\ref{sec:shadow}).
\end{itemize}
That the descent defect tends to $0$ along backward orbits is proved
in Section~\ref{sec:depth}. Suppose
that the descent defect stays above some $\eta > 0$ along a backward orbit.
Moving back along the orbit by near-returns that almost attain the
descent defect, with tolerances that shrink geometrically, we obtain a sequence
$(u_m)$ with $d(u_i, u_j) > \eta/2$ for all $i < j$, which
contradicts (i). Section~\ref{sec:occur} turns to the converse
and asks which defects occur over a given asymptotic distance, and
Section~\ref{sec:summary} summarizes which
properties are forced and which are not.

The theorems, and the stated properties of the three spirals, are
checked in Lean~4 using Mathlib, except for the statements
listed in the Data availability section; version 3 adds
\texttt{Converses.lean}, which checks, among other statements,
Theorem~\ref{thm:pairs} and the maximality of the asymptotic
distance. The Lean files are archived
with this paper and build on the verification files of
\cite{DuzerC}.

\section{Preliminaries}\label{sec:prelim}

A \emph{quasi-pseudo-metric space} is a set $X$ with a function
$d \ge 0$ such that $d(x, x) = 0$ and
$d(x, z) \le d(x, y) + d(y, z)$; neither symmetry nor separation is
assumed. Following \cite{Zav2023}, a \emph{quasi-metric} is a
quasi-pseudo-metric with $d(x, y) > 0$ for $x \ne y$ (the $T_1$
axiom). A quasi-pseudo-metric is \emph{$R_0$} if $d(x, y) = 0$ implies
$d(y, x) = 0$; every quasi-metric is $R_0$. A sequence $(x_n)$ is \emph{left K-Cauchy} if for every
$\varepsilon > 0$ there is $N$ with $d(x_n, x_m) < \varepsilon$
whenever $N \le n \le m$. Throughout, we assume
\begin{itemize}
\item[(i)] every sequence in $X$ has a left K-Cauchy subsequence.
\end{itemize}
For quasi-metric spaces, (i) is equivalent to hereditary
precompactness \cite[Theorem~3]{KMRV1993}, the setting of the title,
but our general results need only (i). For background on asymmetric
spaces, including left K-Cauchy sequences and their topologies, see
the book of Cobza\c{s} \cite{Cobzas2013}.

We also fix a non-contractive bijection $F$ of $X$, with inverse $G$.
Equivalently, $G$ is non-expansive: $d(Gx, Gy) \le d(x, y)$. Iterates
of $F$ are again non-contractive, and iterates of $G$ non-expansive.

A \emph{chain} is a strictly increasing sequence $t_1 < t_2 < \cdots$
of natural numbers; its \emph{differences} are the numbers
$t_j - t_i$ with $i < j$. Hypothesis (i) provides chains along which
an orbit nearly returns to itself, by the mechanism of
\cite[Lemma~3.3]{DuzerC}. Fix $\varepsilon > 0$. Applied to the
forward orbit of a point $b$, (i) gives, after discarding finitely
many terms, a chain with
$d(F^{t_i} b, F^{t_j} b) < \varepsilon$ for $i < j$; since $G^{t_i}$
is non-expansive, $d(b, F^\Delta b) < \varepsilon$ for every
difference $\Delta$. Applied in the same way to the backward orbit of a
point $a$, it gives a chain with $d(G^{t_i} a, G^{t_j} a) < \varepsilon$ for
$i < j$. No upper bound on $d(a, G^\Delta a)$ follows, since $F^{t_i}$ can
only increase distances; this asymmetry is why the source carries
the defect in Section~\ref{sec:source}. What the backward chain does
give is near-returns far back on the orbit:
$d(G^t a, G^{t + \Delta} a) < \varepsilon$ for every $t \ge t_1$ and
arbitrarily large $\Delta$ (Proposition~\ref{prop:deep}).

\begin{definition}[Asymptotic distance, defect, descent defect]\label{def:fee}
For $x, y \in X$, the \emph{ladder} $k \mapsto d(G^k x, G^k y)$ is
nonincreasing and nonnegative. Its limit
\[
\lam(x, y) = \inf_{k \ge 0} d(G^k x, G^k y)
\]
is the \emph{asymptotic distance}, $\D(x, y) = d(x, y) - \lam(x, y) \ge 0$ is
the \emph{defect}, and $\lams(x, y) = \lam(x, y) + \lam(y, x)$ is the
\emph{symmetrized asymptotic distance}. In a pair $(x, y)$, $x$ is the
\emph{source} and $y$ the \emph{target}; the pair is \emph{strict} if
$\D(x, y) > 0$, a \emph{descent} if $y = G^\Delta x$ for some
$\Delta \ge 1$, and an \emph{ascent} if $y = F^\Delta x$ for some
$\Delta \ge 1$. The points $G^\Delta x$ with $\Delta \ge 1$ are the
\emph{backward iterates} of $x$. The \emph{descent
defect}\footnote{The fees of \cite[Section~5]{Duzer2026} are
distances, $d(a, G^\Delta a)$ and $d(G^\Delta b, b)$, not suprema of
defects; we do not use that word here.} of $x$ is the supremum
of the defects of its descents,
\[
c(x) = \sup_{\Delta \ge 1} \D(x, G^\Delta x),
\]
and its \emph{ascent defect} is $A(x) = \sup_{j \ge 1} \D(x, F^j x)$.
\end{definition}

The asymptotic distance is what \cite[Theorem~5.1]{Duzer2026} calls
the invariant limit register, written here in terms of the inverse. We
use the following facts about it, each checked
in the
accompanying file.

\needspace{8\baselineskip}
\begin{lemma}\label{lem:register}
\begin{enumerate}
\item[(a)] $\lam$ is a quasi-pseudo-metric with $\lam \le d$, and
$\lam(Fx, Fy) = \lam(x, y)$.
\item[(b)] $\D(x, y) \le \D(Fx, Fy)$: the defect is nondecreasing
along the dynamics.
\item[(c)] $\D(x, z) \le \D(x, y) + \D(y, z) + \lams(x, y)$.
\item[(d)] Under (i), distances out of a point are bounded, so $c(x)$
and $A(x)$ are finite.
\end{enumerate}
\end{lemma}

\begin{proof}
(a) The triangle inequality passes to the limit along the ladders,
and $\lam \le d$ is the case $k = 0$. For invariance, the ladder of
$(Fx, Fy)$ is that of $(x, y)$ shifted by one step, since
$d(G^{k+1}Fx, G^{k+1}Fy) = d(G^k x, G^k y)$, so the two ladders have
the same limit.
(b) By (a), $\D(Fx, Fy) - \D(x, y) = d(Fx, Fy) - d(x, y) \ge 0$.
(c) Expand $\D$ and use $d(x, z) \le d(x, y) + d(y, z)$ and
$\lam(y, z) \le \lam(y, x) + \lam(x, z)$.
(d) Suppose $d(x, y_k) > k$ for a sequence $(y_k)$. A left K-Cauchy
subsequence $(y_{k_m})$ has an $N$ with $d(y_{k_N}, y_{k_m}) < 1$
for $m \ge N$, so $d(x, y_{k_m}) \le d(x, y_{k_N}) + 1$ stays bounded
while $k_m \to \infty$, a contradiction. Hence $c(x)$ and $A(x)$ are
at most $\sup_y d(x, y) < \infty$.
\end{proof}

\section{The defect is carried by the source}\label{sec:source}

Which endpoint of a strict pair carries its defect? The monochromy
theorem of \cite[Theorem~5.6]{Duzer2026} gives a partial answer. For
one of the two endpoints, it finds a chain along whose differences
$\Delta$ the defect of the pair formed by that endpoint and its
iterate, $\D(a, G^\Delta a)$ for $a$
and $\D(G^\Delta b, b)$ for $b$, is at least a quarter of the pair's
defect. It does not identify the endpoint, and its
proof uses a Ramsey argument. We prove directly that the endpoint is
always the source, and that the quarter can be replaced by any
fraction less than one.

\begin{theorem}[The source carries the defect]\label{thm:source}
Let $a, b \in X$, $\eta \in \R$ and $M \in \N$ with
$d(a, b) - d(G^\Delta a, G^\Delta b) \ge \eta$ for $\Delta \ge M$.
Then for every $\varepsilon > 0$ there is a chain along whose
differences $\Delta$
\[
d(G^\Delta b, b) < \varepsilon
\qquad\text{and}\qquad
d(a, G^\Delta a) > \eta - \varepsilon .
\]
In particular, if $(a, b)$ is strict with defect $\eta_0 = \D(a, b)$,
the hypothesis holds for every $\eta < \eta_0$ and some $M$.
\end{theorem}

In other words, suppose that moving $a$ and $b$ far into the past
brings $a$ at least $\eta$ closer to $b$. Then, for the differences
$\Delta$ of a chain, the distance from $G^\Delta b$ to $b$ is less
than $\varepsilon$, while the distance from $a$ to $G^\Delta a$
exceeds $\eta - \varepsilon$: the target nearly returns to itself,
and the source does not.

\begin{proof}
By (i), the forward orbit $(F^n b)_n$ has a left K-Cauchy
subsequence. Discarding finitely many terms, we get a chain $(t_i)$
with $d(F^{t_i} b, F^{t_j} b) < \varepsilon$ for all $i < j$. For
$\Delta = t_j - t_i$, since $F^{t_j}$ is non-contractive and
$F^{t_j}(G^\Delta b) = F^{t_i} b$,
\[
d(G^\Delta b, b) \le d(F^{t_j} G^\Delta b, F^{t_j} b)
= d(F^{t_i} b, F^{t_j} b) < \varepsilon .
\]
Thin the chain so that all its differences exceed $M$. For such
$\Delta$,
\[
d(a, b) \le d(a, G^\Delta a) + d(G^\Delta a, G^\Delta b)
+ d(G^\Delta b, b)
< d(a, G^\Delta a) + d(G^\Delta a, G^\Delta b) + \varepsilon,
\]
and the hypothesis gives $d(a, G^\Delta a) > \eta - \varepsilon$. For
the last statement, note that $d(G^\Delta a, G^\Delta b)$ decreases
to $\lam(a, b)$.
\end{proof}

\begin{corollary}[The source's orbit is strict and infinite]
\label{cor:orbit}
Let $(a, b)$ satisfy the hypothesis of Theorem~\ref{thm:source} with
$\eta > 0$.
\begin{enumerate}
\item[(a)] For every $\varepsilon > 0$ there are $\Delta \ge 1$ and
$k \ge 0$ with $d(a, G^\Delta a) > \eta - \varepsilon$ and
$d(G^k a, G^{k + \Delta} a) < \varepsilon$. If
$\varepsilon < \eta/2$, the pair $(a, G^\Delta a)$ is strict. In
particular, the orbit of $a$, with the restriction of $d$, is itself
a space satisfying (i) on which $F$ is a non-contractive
non-isometry.
\item[(b)] $G^p a \ne a$ for every $p \ge 1$.
\end{enumerate}
\end{corollary}

\begin{proof}
(a) Take the chain of Theorem~\ref{thm:source} and apply (i) to the
sequence $(G^{t_i} a)_i$: some $i < j$ have
$d(G^{t_i} a, G^{t_j} a) < \varepsilon$, and $\Delta = t_j - t_i$,
$k = t_i$ serve. Then
$\lam(a, G^\Delta a) \le d(G^k a, G^{k+\Delta} a) < \varepsilon$, so
$\D(a, G^\Delta a) > \eta - 2\varepsilon$, which is positive when
$\varepsilon < \eta/2$. Hypothesis (i) passes to subspaces, $F$ is
non-contractive on the orbit because it is on $X$, and the strict
pair $(a, G^\Delta a)$ lies in the orbit, so $F$ is not an isometry
there.
(b) If $G^p a = a$, take the chain of Theorem~\ref{thm:source} with
$\varepsilon = \eta/2$. Its terms take only finitely many values
modulo $p$, so $t_i \equiv t_j \pmod p$ for some $i < j$. Then
$\Delta = t_j - t_i$ is a difference of the chain and a multiple of
$p$, so $G^\Delta a = a$ and $0 = d(a, G^\Delta a) > \eta/2$, a
contradiction.
\end{proof}

\begin{corollary}[The source carries almost all of the defect]\label{cor:uniform}
Let $(a, b)$ satisfy the hypothesis of Theorem~\ref{thm:source}.
Then for every $\varepsilon > 0$ there is a chain along whose
differences $\Delta$
\[
\D(a, G^\Delta a) > \eta - 2\varepsilon .
\]
In particular, if $(a, b)$ is strict, then for every
$\varepsilon' > 0$ there is a chain along whose differences
$\D(a, G^\Delta a) > \D(a, b) - \varepsilon'$.
\end{corollary}

\begin{proof}
Take the chain $(t_i)$ of Theorem~\ref{thm:source} and refine it: by
(i) applied to the sequence $(G^{t_i} a)_i$, there is a subchain
along which $d(G^{t_i} a, G^{t_j} a) < \varepsilon$ for $i < j$. For
a difference $\Delta = t_j - t_i$ of the subchain,
\[
\lam(a, G^\Delta a) \le d\bigl(G^{t_i} a, G^{t_i}(G^\Delta a)\bigr)
= d(G^{t_i} a, G^{t_j} a) < \varepsilon,
\]
while $d(a, G^\Delta a) > \eta - \varepsilon$ by
Theorem~\ref{thm:source}. Subtracting gives the claim.
\end{proof}

Corollary~\ref{cor:uniform} improves the monochromy theorem of
\cite{Duzer2026} in both respects: the endpoint is the source, and
the bound is arbitrarily close to the full defect.

Corollary~\ref{cor:orbit} also gives a normal form. Call the orbit of
the source of a strict pair a \emph{strict orbit}. By (b), its points
$F^n a$, $n \in \Z$, are distinct, so $n \mapsto F^n a$ identifies it
with $\Z$ and $F$ with the shift $n \mapsto n + 1$. When $d$ is a
quasi-metric, its restriction is a hereditarily precompact
quasi-metric on $\Z$ for which the shift is non-contractive but not
an isometry. The normal form applies orbit by orbit, not to $F$ as a
whole. The spirals of Section~\ref{sec:spirals} and the zone tolls of
Section~\ref{sec:depth} are written in this form.

\needspace{10\baselineskip}
\begin{proposition}[Upper bound by the descent defect]\label{prop:fee}
\begin{enumerate}
\item[(a)] For all $a, b$ and $\varepsilon > 0$ there is $\Delta \ge 1$
with $\D(a, G^\Delta a) > \D(a, b) - \varepsilon$.
\item[(b)] $\D(a, b) \le c(a)$ for all $a, b$.
\item[(c)] $c(a) \le c(Fa)$: the descent defect is nondecreasing along orbits.
\end{enumerate}
\end{proposition}

\begin{proof}
(a) If $\D(a, b) < \varepsilon$, any $\Delta$ serves, since defects
are nonnegative. Otherwise, since the ladder of $(a, b)$ is
nonincreasing, there is $K$ with
$d(G^k a, G^k b) < \lam(a, b) + \varepsilon/3$ for $k \ge K$. The
hypothesis of Theorem~\ref{thm:source} then holds with
$\eta = \D(a, b) - \varepsilon/3 > 0$ and $M = K$, and
Corollary~\ref{cor:orbit}(a), applied with $\varepsilon/3$, gives
$\Delta$ and $k$ with $d(a, G^\Delta a) > \eta - \varepsilon/3$ and
$\lam(a, G^\Delta a) \le d(G^k a, G^{k+\Delta} a) < \varepsilon/3$.
Subtracting gives
$\D(a, G^\Delta a) > \eta - 2\varepsilon/3 = \D(a, b) - \varepsilon$.
(b) Take the supremum over $\Delta$ in (a).
(c) By Lemma~\ref{lem:register}(b),
$\D(a, G^\Delta a) \le \D(Fa, F G^\Delta a) = \D(Fa, G^\Delta Fa)$;
take the supremum over $\Delta$.
\end{proof}

Proposition~\ref{prop:fee} gives the upper bound in
Theorem~\ref{thm:structure}(b), together with the monotonicity of
the descent defect: every defect is bounded by the descent defect of its source, which
does not decrease in time. For example, the toll spiral, whose toll
is a constant $\kappa > 0$ on descents from
points $n \ge 0$, has descent defect $\kappa$ at $n \ge 0$ and
$0$ at $n < 0$, and the height spiral of Section~\ref{sec:spirals}
has descent defect $2^n$ at $n$.

\section{Upper and lower bounds for the defect}\label{sec:band}

Proposition~\ref{prop:fee}(b) bounds the defect of a pair above by
the descent defect of its source. The matching lower bound involves the descent defects of
both endpoints, with the symmetrized asymptotic distance as a correction term, and it holds
not only now but at every earlier time.

\begin{theorem}[Two-sided bound for the defect]\label{thm:band}
For all $x, y \in X$ and $k \ge 0$,
\[
c(G^k x) - c(G^k y) - \lams(x, y) \;\le\; \D(x, y) \;\le\; c(x) .
\]
\end{theorem}

\begin{proof}
The upper bound is Proposition~\ref{prop:fee}(b). For the lower bound
with $k = 0$, let $\varepsilon > 0$ and choose $\Delta \ge 1$ with
$\D(x, G^\Delta x) > c(x) - \varepsilon$. Lemma~\ref{lem:register}(c),
and Proposition~\ref{prop:fee}(b) applied to the pair
$(y, G^\Delta x)$, give
\[
c(x) - \varepsilon < \D(x, G^\Delta x)
\le \D(x, y) + \D(y, G^\Delta x) + \lams(x, y)
\le \D(x, y) + c(y) + \lams(x, y).
\]
For general $k$, apply this to the pair $(G^k x, G^k y)$. Its defect
is at most $\D(x, y)$ by Lemma~\ref{lem:register}(b), applied $k$
times, and its symmetrized asymptotic distance equals $\lams(x, y)$ by
invariance.
\end{proof}

Both bounds are attained, by the toll and thin spirals of
Section~\ref{sec:spirals} (see the paragraph after
Proposition~\ref{prop:thin}). Since $\D \ge 0$, the lower bounds for
all $k$ together say that
\[
\D(x, y) \;\ge\; \bigl(\tau_c(x, y) - \lams(x, y)\bigr)^+,
\qquad
\tau_c(x, y) = \sup_{k \ge 0}\,\bigl(c(G^k x) - c(G^k y)\bigr)^+ .
\]
The function $\tau_c$ depends only on $c$. By
Theorem~\ref{thm:occur}, it is itself the defect of $F$ for the
distance $\lam + \tau_c$, which has the same asymptotic distance and
the same descent defect.

The next proposition compares the descent defect with its
counterpart for ascents, the ascent defect
$A(x) = \sup_{j \ge 1} \D(x, F^j x)$ of Definition~\ref{def:fee}: the
ascent defect never exceeds the descent defect, and the two have the
same limit in the deep past.

\needspace{10\baselineskip}
\begin{proposition}[Descent and ascent defects]\label{prop:afee}
\begin{enumerate}
\item[(a)] $A(x) \le c(x)$ for every $x$, and $A(x) \le A(Fx)$.
\item[(b)] For every $a \in X$ and $\varepsilon > 0$ there is a chain
$(t_r)$ such that for every $r$,
\[
c(G^{t_r} a) \le A(G^{t_{r+1}} a) + 2\varepsilon .
\]
\item[(c)] The backward limits $c_\infty(a) = \lim_k c(G^k a)$ and
$A_\infty(a) = \lim_k A(G^k a)$ exist, since $c$ and $A$ are
nondecreasing along orbits, and they coincide. (By
Theorem~\ref{thm:vanish}, both are $0$.)
\end{enumerate}
\end{proposition}

\begin{proof}
(a) The defect $\D(x, F^j x)$ of each ascent is at most $c(x)$ by
Proposition~\ref{prop:fee}(b), and $A(x) \le A(Fx)$ follows from
Lemma~\ref{lem:register}(b).

(b) Apply (i) first to the backward orbit of $a$, and then, along the
indices obtained, to the forward orbit. After discarding finitely
many terms, this gives a chain $(t_r)$ with
$d(G^{t_r} a, G^{t_s} a) < \varepsilon$ and
$d(F^{t_r} a, F^{t_s} a) < \varepsilon$ for all $r < s$.

Fix $r < s$ and write $p = t_r$, $q = t_s$. Since $G^{p+q}$ maps
$(F^p a, F^q a)$ to $(G^q a, G^p a)$, invariance of the asymptotic distance
gives
\[
\lam(G^q a, G^p a) = \lam(F^p a, F^q a) \le d(F^p a, F^q a)
< \varepsilon .
\]
For any $z$, combining this with the triangle inequalities
$d(G^p a, z) \le d(G^p a, G^q a) + d(G^q a, z)$ and
$\lam(G^q a, z) \le \lam(G^q a, G^p a) + \lam(G^p a, z)$ gives
\[
\D(G^p a, z) < \D(G^q a, z) + 2\varepsilon .
\]

Now let $z = G^{\Delta}(G^{t_r} a)$ be a backward iterate of
$G^{t_r} a$, and choose $s > r + 1$ with $t_s > t_r + \Delta$. Then
$z$ is a forward iterate of $G^{t_s} a$, so
$\D(G^{t_s} a, z) \le A(G^{t_s} a) \le A(G^{t_{r+1}} a)$ by (a), and
the last display gives
$\D(G^{t_r} a, z) < A(G^{t_{r+1}} a) + 2\varepsilon$. Taking the
supremum over $\Delta$ proves (b).

(c) By (a), $A_\infty \le c_\infty$. By (b),
$c_\infty(a) \le c(G^{t_r} a) \le A(G^{t_{r+1}} a) + 2\varepsilon$
for every $r$, and letting $r \to \infty$ gives
$c_\infty \le A_\infty + 2\varepsilon$.
\end{proof}

Finally, given $\varepsilon > 0$, every point far enough back on an
orbit has descents of arbitrarily large displacement and of cost
less than $\varepsilon$. Hypothesis (i) provides one such
point on each orbit, and the property passes to all earlier points
because $G$ is non-expansive.

\begin{proposition}[Near-returns far back on an orbit]\label{prop:deep}
For every $a \in X$ and $\varepsilon > 0$ there is $T$ such that for
every $t \ge T$ and every $N$ there is $\Delta \ge N$ with
$d(G^t a, G^{t + \Delta} a) < \varepsilon$.
\end{proposition}

\begin{proof}
By Section~\ref{sec:prelim}, the backward orbit of $a$ has a chain
$(t_i)$ with $d(G^{t_i} a, G^{t_j} a) < \varepsilon$ for $i < j$. Let
$T = t_1$. Given $t \ge T$ and $N$, pick $j$ with $t_j - t_1 \ge N$
and put $\Delta = t_j - t_1$. Since $G^t a = G^{t - t_1}(G^{t_1} a)$,
$G^{t+\Delta} a = G^{t - t_1}(G^{t_j} a)$ and $G^{t - t_1}$ is
non-expansive,
$d(G^t a, G^{t + \Delta} a) \le d(G^{t_1} a, G^{t_j} a) < \varepsilon$.
\end{proof}

\section{The asymptotic distance on an orbit}\label{sec:shadow}

On the orbit of a point $a$, the asymptotic distance depends only on the
displacement: by invariance, $\lam(F^n a, F^m a)$ is a function of
$m - n$. Since $\lam \le d$, it also satisfies (i). These
two facts are enough to show that the symmetrized asymptotic distance is
totally bounded on each orbit. We prove this in the
following sequential form, from which total boundedness follows in
the usual way.

\begin{theorem}[The symmetrized asymptotic distance is totally bounded on orbits]\label{thm:shadow}
Let $a \in X$ and write $a_n = F^n a$ for $n \in \Z$, with $F^{-k} =
G^k$. For every sequence $(n_i)$ of integers and every $\varepsilon >
0$ there are $i < j$ with
\[
\lam(a_{n_i}, a_{n_j}) + \lam(a_{n_j}, a_{n_i}) < \varepsilon .
\]
\end{theorem}

\begin{proof}
Put $\ell(n, m) = \lam(a_n, a_m)$. By Lemma~\ref{lem:register}(a),
$\ell(n + 1, m + 1) = \ell(n, m)$, so $\ell(n, m) = f(m - n)$ for a
function $f$ on $\Z$; and $\ell$ satisfies (i) because $\ell \le d$.
Let $(n_i)$ be a sequence. By (i) for $\ell$, as in
Section~\ref{sec:prelim}, there is a chain $(\mu_i)$ with
$\ell(n_{\mu_i}, n_{\mu_j}) < \varepsilon/2$ for $i < j$.
Applying (i) for $\ell$ to the sequence $(-n_{\mu_i})_i$ gives a
further chain $(\nu_i)$ with
$\ell(-n_{\mu_{\nu_i}}, -n_{\mu_{\nu_j}}) < \varepsilon/2$
for $i < j$. Since $\ell(-u, -v) = f(u - v) = \ell(v, u)$, the
integers $n_{\mu_{\nu_i}}$ and $n_{\mu_{\nu_j}}$, $i < j$,
are within $\varepsilon/2$ of each other in both directions.
\end{proof}

We now complete the orbit. In the notation of the proof, $f$ is
nonnegative and subadditive with $f(0) = 0$, so
$f_s(k) = f(k) + f(-k)$ is a symmetric subadditive function, and
$(n, m) \mapsto f_s(m - n)$ is an invariant pseudo-metric on $\Z$ (the
symmetrized asymptotic distance on the orbit). By Theorem~\ref{thm:shadow}, it
is totally bounded. Its kernel $\{k : f_s(k) = 0\}$ is a subgroup
$p\Z$. By Weil's theorem on precompact abelian groups
\cite[\S 3.7]{AT2008}, completing $\Z/p\Z$ with respect to $f_s$ yields a compact
metrizable abelian group $K$, finite if $p > 0$. The image of $1$
generates a dense subgroup of $K$, so $K$ is monothetic, and the
shift acts on it as the rotation by this generator.

Since $|f(u) - f(v)| \le f_s(u - v)$, the function $f$ extends
continuously to $K$, and $\ell$ extends to an invariant
quasi-pseudo-metric $\bar\ell$ on $K$ whose symmetrization is the
metric of $K$. A priori, the topology of $\bar\ell$ is only known to
be coarser than the compact topology, so $K$ is merely a totally
bounded abelian paratopological group in which the image of
$(\Z, \ell)$ is dense. Completions of this kind are studied in
\cite{KRS1998}, and totally bounded paratopological groups in
\cite{Tk2013}. We call $(K, \bar\ell)$, with its rotation, the
\emph{completed orbit} of $a$. Apart from Corollary~\ref{cor:tollanswer}, no result
outside this section uses the completion.
The next theorem shows that its topology is the compact one.

\begin{theorem}[The completed orbit is a topological group]\label{thm:shadtop}
The topology of $\bar\ell$ equals the compact topology of $K$. In
particular, the completion is a compact topological group, not merely a
paratopological one.
\end{theorem}

\begin{proof}
The extension of $f$ to $K$ is continuous (indeed $1$-Lipschitz for
the metric $f_s$), subadditive and nonnegative, with $f(0) = 0$. Hence
$H = \{x \in K : f(x) = 0\}$ is a closed subsemigroup of $K$.

A closed subsemigroup of a compact group is a subgroup. Indeed, let
$x \in H$. By compactness, the multiples $nx$, $n \ge 1$, have a
convergent subsequence $(n_k x)$, and the differences
$m_k = n_{k+1} - n_k$, which we may take $\ge 2$, satisfy
$m_k x \to 0$. Then $(m_k - 1)x \in H$ and $(m_k - 1)x \to -x$, so
$-x \in H$ because $H$ is closed. Therefore
$f_s(x) = f(x) + f(-x) = 0$, and $x = 0$ since $f_s$ is a metric on $K$.
Thus $H = \{0\}$.

The closed sets $\{f \le 1/n\}$ decrease and have intersection
$H = \{0\}$, so by compactness they eventually lie inside any
$K$-neighborhood of $0$. Conversely, $f$ is continuous, so every
$\bar\ell$-ball is $K$-open. Hence the two topologies have the same
neighborhoods of $0$. Both are invariant under translation, so they
coincide.
\end{proof}

Conversely, every such group, quasi-pseudo-metric and rotation
arises in this way, even from a map that is not an isometry.

\begin{proposition}[Every completion occurs]\label{prop:realize}
Let $K$ be a compact metrizable abelian group, $g \in K$ an element
generating a dense subgroup, and $\ell$ an invariant
quasi-pseudo-metric on $K$ whose symmetrization $\ell_s$ is a metric
inducing the topology of $K$. For $\kappa > 0$, define on $\Z$
\[
d(n, m) = \ell(ng, mg) + \tau(n, m), \qquad
\tau(n, m) =
\begin{cases}
\kappa & \text{if } m < n \text{ and } n \ge 0,\\
0 & \text{otherwise.}
\end{cases}
\]
Then $(\Z, d)$ is a quasi-pseudo-metric space satisfying (i), and a
$T_1$ quasi-metric space if $K$ is infinite. The shift $S(n) = n + 1$ is a
non-contractive bijection of it that is not an isometry. Its asymptotic distance
is $\lam(n, m) = \ell(ng, mg)$, and the completed orbit of $0$ is
$(K, \ell)$, where the class of $n$ is sent
to $ng$ and $S$ acts as the rotation by $g$.
\end{proposition}

\begin{proof}
The base distance $(n, m) \mapsto \ell(ng, mg)$ is a quasi-pseudo-metric, and
it is invariant under the shift. For the toll, suppose that
$\tau(n, m) = \kappa$, so that $m < n$ and $n \ge 0$, and let
$k \in \Z$. If $k < n$, then $\tau(n, k) = \kappa$. Otherwise
$k \ge n \ge 0$ and $m < k$, so $\tau(k, m) = \kappa$. Hence
$\tau(n, m) \le \tau(n, k) + \tau(k, m)$, and $d$ is a
quasi-pseudo-metric. Moreover, $\ell(x, y) > 0$ for $x \ne y$. Indeed,
$f(x) = \ell(0, x)$ is continuous on $K$, by the inequality at the end
of this proof, and it is subadditive and nonnegative with
$f(0) = 0$. As in the proof of Theorem~\ref{thm:shadtop},
$\{f = 0\}$ is then a closed subgroup of $K$, and $\ell_s(0, x) = 0$
on it, so it is $\{0\}$. If $K$ is infinite, then $g$ has infinite
order, since a finite subgroup is closed and could not be dense. So
$d(n, m) \ge \ell(ng, mg) > 0$ for $n \ne m$, and $d$ is a $T_1$
quasi-metric.

The shift maps descents to descents and nonnegative integers to
nonnegative integers, so $\tau(n + 1, m + 1) \ge \tau(n, m)$, and $S$
is non-contractive. It is not an isometry, since
$d(0, -1) = d(-1, -2) + \kappa$.

For (i), let $(n_i)$ be a sequence in $\Z$. If some value occurs
infinitely often, it gives a constant subsequence. Otherwise, since
$K$ is compact, there is a strictly monotone subsequence along which
$n_i g$ converges in $K$, and along it
$\ell(n_i g, n_j g) \le \ell_s(n_i g, n_j g) \to 0$. If the
subsequence increases, each pair $(n_i, n_j)$ with $i < j$ is an
ascent and pays no toll. If it decreases, its terms are eventually
negative, and a pair whose source is negative pays no toll. In both
cases the subsequence is left K-Cauchy.

For $k > n$, the integer $n - k$ is negative, so
$d(n - k, m - k) = \ell(ng, mg)$ by invariance of $\ell$. Hence
$\lam(n, m) = \ell(ng, mg)$, and the defect is the toll. In the
notation of this section, on the orbit of $0$ we have
$f(k) = \ell(0, kg)$ and $f_s(k) = \ell_s(0, kg)$, so the kernel $p\Z$
is $\{k : kg = 0\}$, and $n + p\Z \mapsto ng$ is an injective
homomorphism of $\Z/p\Z$ into $K$ and an isometry from the metric
$f_s(m - n)$ to $\ell_s$. Its
image is dense, and $K$ is complete, being a compact metric space.
So the completion is $K$, the image of $1$ is $g$, and the shift acts
as the rotation by $g$. Finally, $\bar\ell = \ell$, because $\ell$ is
continuous on $K$:
$|\ell(x, y) - \ell(x', y')| \le \ell_s(x, x') + \ell_s(y, y')$.
\end{proof}

With the construction of this section, the proposition shows that,
up to isomorphism, the completed orbits are exactly the compact metrizable
monothetic groups with an invariant quasi-pseudo-metric whose
symmetrization is a metric inducing their topology, together with the
rotation by a generator. Every compact metrizable abelian group has
such a quasi-pseudo-metric, and even a metric: if $e$ is any metric
inducing its topology, then $\sup_{z \in K} e(x + z, y + z)$ is an
invariant metric, finite because $e$ is bounded, and inducing the same
topology because $(w, v) \mapsto e(v + w, v)$ is uniformly continuous
on $K \times K$ and vanishes for $w = 0$. For
$K = \mathbb{T}$ with the arc metric and $g = \theta$ irrational, the
space of the proposition is the toll spiral of \cite{Duzer2026},
described in Section~\ref{sec:spirals}.

\begin{remark}[An asymmetric example]\label{rem:asym}
The asymptotic distance need not be a metric. On the circle, let
$\varphi(x, y) = \min(t, 2(1 - t))$, where $t$ is the length of the
counter-clockwise arc from $x$ to $y$. This is an invariant $T_1$
quasi-metric, and its symmetrization is $\min(3t, 3(1 - t), 1)$, a
metric inducing the usual topology. By
Proposition~\ref{prop:realize} with $g = \theta$ irrational,
$(n, m) \mapsto \varphi(n\theta, m\theta)$ is the asymptotic distance of a
non-contractive bijection of a $T_1$ quasi-metric space satisfying
(i), and the completed orbit of $0$ is the circle with $\varphi$. So the
asymptotic distance can be asymmetric, although by
Theorem~\ref{thm:shadtop} the topology of the completed orbit is always
the compact group topology. What Theorem~\ref{thm:shadow} does
rule out, on an infinite orbit, is the one-way circle, in which $\ell(x, y)$
is the length of the counter-clockwise arc from $x$ to $y$:
there $\ell(x, y) + \ell(y, x) = 1$ for $x \ne y$, so the symmetrization is
discrete.
\end{remark}

\section{Three examples on \texorpdfstring{$\Z$}{Z}, and tolls over a shadow}\label{sec:spirals}

This section concerns \cite[Question~7.6]{Duzer2026}: ``Is every
bijective non-contractive map that is not an isometry, on a
hereditarily precompact quasi-metric space, conjugate in an
appropriate sense to a tolled extension (a fixed toll zone over an
isometric shadow), possibly with infinitely many zones?''
Sections~\ref{sec:height}--\ref{sec:ascent} give three examples,
Section~\ref{sec:tolled} characterizes tolled extensions on an
orbit, and Section~\ref{sec:tollshadow} gives the structural
answer.

Throughout, $\mathbb{T} = \R/\Z$ carries the arc metric, $\theta$ is
irrational, $\varphi(n) = n\theta$, and $\rho(n, m)$ denotes the arc
distance between $\varphi(n)$ and $\varphi(m)$. Then $\rho$ is
shift-invariant, $\rho(n + 1, m + 1) = \rho(n, m)$, and vanishes only
on the diagonal. Each spiral below is a quasi-metric $d = \rho + \tau$
on $\Z$, where the \emph{toll} $\tau \ge 0$ satisfies
$\tau(n + 1, m + 1) \ge \tau(n, m)$, so that the shift
$S(n) = n + 1$ is non-contractive. In the notation of
Section~\ref{sec:prelim}, a pair $(n, m)$ is a descent if $m < n$ and
an ascent if $m > n$. We call $n$ \emph{recent} if $n \ge 0$ and
\emph{deep} if $n < 0$. In all three spirals the toll of every pair
tends to $0$ as the pair is moved into the past, so the asymptotic distance is
$\rho$ and the defect is the toll.

In this notation, the toll spiral of \cite{Duzer2026} has
$\tau(n, m) = \kappa$, a constant $\kappa > 0$ (written $c$ in
\cite{Duzer2026}), on descents from recent sources ($m < n$ and
$n \ge 0$) and $\tau = 0$ otherwise. In the terms of the question,
its isometric shadow is $(\Z, \rho)$ and its toll zone is the set of
descents from recent sources. Its toll depends only on the
source, is subadditive on its own, and its descent defect is the
step function $c(n) = \kappa$ for $n \ge 0$ and $c(n) = 0$ for $n < 0$.

\subsection{The height spiral: the toll can depend on the target}\label{sec:height}

\begin{example}[Height spiral]\label{ex:height}
Let $h(n) = 2^n$ and
\[
d(n, m) = \rho(n, m) + \bigl(h(n) - h(m)\bigr)^+ .
\]
\end{example}

\begin{proposition}\label{prop:height}
The height spiral is a $T_1$ quasi-metric on $\Z$ satisfying (i). The
shift is non-contractive and strictly expands every descent. The toll
depends on the target: $\tau(1, 0) \ne \tau(1, -1)$. The descent defect is $c(n)
= 2^n$.
\end{proposition}

\begin{proof}
The toll $(h(n) - h(m))^+$ is subadditive, since
$h(n) - h(k) = (h(n) - h(m)) + (h(m) - h(k))$ and the positive part
is subadditive. So $d$ is a quasi-pseudo-metric, and it is $T_1$
because $\rho$ separates points. The shift doubles the toll, so it is
non-contractive, and strictly expanding on every descent, where the
toll is positive. Also $\tau(1, 0) = 1$ and $\tau(1, -1) = 3/2$.

For (i), take a sequence $(u_k)$ and pass to a subsequence along which
$\varphi(u_k)$ converges and $(u_k)$ is monotone. If $(u_k)$ is
nondecreasing, the toll from each term to any later one vanishes, and
the base is Cauchy. If it is strictly decreasing, then
$u_k \to -\infty$, and the toll from $u_k$ to any later term is at
most $h(u_k) \to 0$.

Finally, $(h(n - k) - h(m - k))^+ \to 0$ as $k \to \infty$, so the
asymptotic distance is $\rho$. Hence $\D(n, n - \Delta) = h(n) - h(n - \Delta)$,
and $c(n) = 2^n$.
\end{proof}

\begin{remark}[Readings of {\cite[Question~7.6]{Duzer2026}}]\label{rem:readings}
The question allows infinitely many zones but does not say how such
a toll is charged. In each natural reading only descents are
charged: the toll of a descent $(n, m)$ is a function $\alpha(n)$ of its
source (the \emph{source-indexed} reading), a sum
$\sum_{m < j \le n} \gamma_j$ of charges $\gamma_j$ at the points $j$
it passes over (the \emph{boundary-crossing} reading), or the sum of the
two (the \emph{mixed} reading), with $\alpha(n), \gamma_j \ge 0$. The height
toll is a boundary-crossing toll, since
\[
(2^n - 2^m)^+ = \sum_{m < j \le n} 2^{j-1},
\]
so it refutes only the source-indexed reading: the descents from $1$
to $0$ and from $1$ to $-1$ are charged $1$ and $3/2$, which cannot
both equal $\alpha(1)$.

The thin and ascent spirals (Examples~\ref{ex:thin}
and~\ref{ex:ascent} below) refute all three readings. Suppose that a
bijection $\sigma$ of $\Z$ commuting with the shift carries the
distance $d$ of a spiral onto a distance $d' = \ell + \tau$, that is,
$d'(\sigma n, \sigma m) = d(n, m)$, with $\ell$ a shift-invariant
quasi-pseudo-metric and $\tau$ a mixed toll. Then the shift is
non-contractive for $d'$, and $d'$ inherits (i), since $\sigma$
carries left K-Cauchy sequences to left K-Cauchy sequences. By (i), $\alpha(n) \to 0$ as $n \to -\infty$: otherwise some
descending sequence has $\alpha(n) \ge \eta > 0$, so each of its terms
is at distance at least $\eta$ from every later one, and no
subsequence is left K-Cauchy. For the same reason, the charges
$\gamma_j$ are summable over $j \le 0$. The toll of the descent
$(n - k, m - k)$ therefore tends to $0$ as $k \to \infty$, which
means that the asymptotic distance of $d'$ is $\ell$ and that its defect is the
toll $\tau$. Hence ascents carry no defect, and moving the target
of a descent deeper only adds charges $\gamma_j$, so its defect does
not decrease. Now $\sigma$ preserves defects: it carries $d$ onto
$d'$ and commutes with the shift, so it carries ladders to ladders.
Being a bijection of $\Z$ that commutes with the shift, it is also a
translation, so it preserves the order of the integers, and with it
ascents and the order of targets. In the thin spiral, however, the
defect $(\kappa - 2\rho)^+$ of a descent from a recent source $n$
can decrease as the target moves deeper: if $\rho(n, m) < \kappa/2$,
then, since the points $\varphi(m')$ with $m' < m$ are dense, there
is $m' < m$ with $\rho(n, m') > \rho(n, m)$, and the descent to $m'$
has smaller defect. So no such
$\sigma$ carries the thin spiral onto a mixed toll, and none carries
the ascent spiral, whose ascents carry defect. Under (i) and with a
non-contractive shift, a toll of any of the three readings gives a
tolled extension with descent zones, and the thin spiral is one too
(Proposition~\ref{prop:tolled} below). The ascent spiral, by contrast,
lies outside this larger class, even up to bi-Lipschitz conjugacy (Proposition~\ref{prop:bilip} below).

The height and thin spirals do not refute the mixed reading up to bi-Lipschitz
conjugacy: the height toll is itself a mixed toll, and since the thin
spiral has $d \ge \max(\rho, \kappa/2)$ on descents from recent sources,
while the toll spiral has $d = \rho + \kappa$ there, the identity is a
$3$-bi-Lipschitz map from the thin spiral onto the toll spiral.
\end{remark}

The height spiral also shows that a toll graded by a quantity
monotone in time can produce a non-isometric example, although every such design tried in
\cite[Remark~5.8]{Duzer2026} failed. Its descent defect $2^n$ tends
to $0$ along backward orbits, its shift strictly expands every
descent, its asymptotic distance $\rho$ separates points, and its toll is
subadditive. The integer layer of the space of
\cite[Proposition~6.1]{DuzerC}, on which a descent from $n$ is
charged $g(n)$ for a positive nondecreasing $g$ tending to $0$ in
the deep past, is another example.

Unlike the toll and thin spirals, the height spiral has an unbounded
toll, and this is forced for every nonzero toll of the form
$(h(n) - h(m))^+$ with $h$ nondecreasing. Indeed, non-contractivity
of the shift then makes the increments of $h$ nondecreasing, and a
nonconstant nondecreasing $h$ with nondecreasing increments is
unbounded. The same holds for the ascent
spiral, whose toll on descents has this form.

\subsection{The thin spiral: a source with positive descent defect can have descents without defect}\label{sec:thin}

\begin{example}[Thin spiral]\label{ex:thin}
Let $\kappa > 0$ and
\[
d(n, m) =
\begin{cases}
\max\bigl(\rho(n, m),\, \kappa - \rho(n, m)\bigr) & \text{if } m < n
\text{ and } n \ge 0,\\
\rho(n, m) & \text{otherwise.}
\end{cases}
\]
Equivalently, $\tau(n, m) = (\kappa - 2\rho(n, m))^+$ on descents from
recent sources and $\tau = 0$ elsewhere.
\end{example}

\begin{proposition}\label{prop:thin}
The thin spiral is a $T_1$ quasi-metric on $\Z$ satisfying (i), and
the shift is non-contractive. Moreover:
\begin{enumerate}
\item[(a)] There are $j \ge 1$ with $\rho(0, -j) < \kappa/2$. For each such
$j$, the descent from $0$ to $-j$ costs more than its base distance,
and the shift strictly expands the pair $(-1, -1 - j)$; in
particular, the shift is not an isometry.
\item[(b)] If $m < 0$ and $\rho(0, m) \ge \kappa/2$, the descent from $0$
to $m$ costs exactly its base distance.
\item[(c)] The descent defect is the same step function as for the
toll spiral: $c(n) = \kappa$ for $n \ge 0$ and $c(n) = 0$ for
$n < 0$. Every descent $(x, y)$ from a recent point $x$ has defect
$(\kappa - 2\rho(x, y))^+$.
\end{enumerate}
\end{proposition}

Since $\theta$ is irrational, $\rho < 1/2$ between distinct points, and the
backward orbit of $0$ is dense, targets as in (b) exist exactly when
$\kappa < 1$.

\begin{proof}
Write $\rho_{nk}$ for $\rho(n, k)$. Since the toll is not
subadditive, the triangle inequality for $d$ has to use the base
distance as well. Let $n \ge 0$ and $k < n$, so that
$d(n, k) = \max(\rho_{nk}, \kappa - \rho_{nk})$, and let $m$ be any
intermediate point. Since $d \ge \rho$, it suffices to show
$\kappa - \rho_{nk} \le d(n, m) + d(m, k)$.
\begin{itemize}
\item If $m < n$, the first leg is a descent from a recent source, so
$d(n, m) \ge \kappa - \rho_{nm}$, and
$\kappa - \rho_{nk} \le (\kappa - \rho_{nm}) + \rho_{mk}$.
\item If $m \ge n$, then $m \ge 0$ and $k < m$, so the second leg is
a descent from a recent source, $d(m, k) \ge \kappa - \rho_{mk}$, and
$\kappa - \rho_{nk} \le \rho_{nm} + (\kappa - \rho_{mk})$.
\end{itemize}
In every other triangle, the left-hand side is a base distance, and
the inequality holds because $d \ge \rho$. Also $d$ is $T_1$ because
$d \ge \rho$, and the shift is non-contractive because a descent from
a recent source remains one. For (i), pass as before to a subsequence
with convergent base and monotone indices: nondecreasing indices
carry no toll, and decreasing indices tend to $-\infty$, after which
the sources are deep and their descents carry no toll.

Since $\theta$ is irrational, backward iterates of $0$ come
arbitrarily close to $0$ in the base, which gives the $j$. For such
$j$, $d(0, -j) = \kappa - \rho_{0,-j} > \rho_{0,-j}$, while
$d(-1, -1 - j) = \rho_{0,-j}$ by invariance. If $\rho(0, m) \ge \kappa/2$,
the maximum is the base. The asymptotic distance is $\rho$, because descents
from deep sources carry no toll. Hence $\D(x, y) = (\kappa - 2\rho(x, y))^+$ for every
descent $(x, y)$ from a recent source $x$, and the descent defect $c(n)$, the supremum of
the defects of the descents from $n$, is $\kappa$ for recent $n$
(approached along close returns) and $0$ for deep $n$.
\end{proof}

So the thin spiral and the toll spiral have the same descent defect
and the same asymptotic distance, and on descents from recent to deep points
between them they attain both bounds of Theorem~\ref{thm:band}. For such a
pair $(x, y)$, the lower bound is $\kappa - 2\rho(x, y)$ (since
$c(y) = 0$ and $\lams(x, y) = 2\rho(x, y)$) and the upper bound is
$c(x) = \kappa$. The toll spiral has defect exactly $\kappa$, the
upper bound, and the thin spiral has defect
$(\kappa - 2\rho(x, y))^+$, the lower bound when
$\rho(x, y) \le \kappa/2$. A source with positive descent defect may
still have descents with zero defect. Take, for example,
$\kappa \le 2\rho(0, -1)$. Then the source $0$ has descent defect
$\kappa$, but its descent to $-1$ has zero defect: it costs exactly
the base distance. The instance checked in
Lean is of this
kind.

\begin{proposition}[The thin toll is the smallest of its kind]
\label{prop:extremal}
Let $d = \rho + \tau$ be a quasi-pseudo-metric on $\Z$, with
$\tau \ge 0$, such that descents from deep sources carry no toll. Let
$n$ be recent, and suppose that there is a nonincreasing function
$\tau_n : (0, 1/2) \to [0, \infty)$ with $\tau(n, m) = \tau_n(\rho(n, m))$
for every deep $m$, and $\tau_n(r) \to \kappa$ as $r \to 0$. Then
$\tau_n(r) \ge (\kappa - 2r)^+$ for every $r \in (0, 1/2)$.
\end{proposition}

\begin{proof}
Since $\tau_n \ge 0$, it suffices to show $\tau_n(r) \ge \kappa - 2r$.
Let $r \in (0, 1/2)$ and $\delta > 0$. Since $\theta$ is irrational,
the points $\varphi(m)$ with $m$ below any given integer are dense in
$\mathbb{T}$. Choose a deep $m$ with $r \le \rho(n, m) < r + \delta$,
and then $k < m$ with $\rho(n, k) < \delta$. Both $m$ and $k$ are
deep, so the descent $(m, k)$ carries no toll, and the triangle
inequality for the route from $n$ through $m$ to $k$ gives
\[
\rho(n, k) + \tau_n(\rho(n, k)) \le \rho(n, m) + \tau_n(\rho(n, m))
+ \rho(m, k).
\]
Here $\rho(n, m) < r + \delta$ and
$\rho(m, k) \le \rho(m, n) + \rho(n, k) < r + 2\delta$, and, since
$\tau_n$ is nonincreasing, $\tau_n(\rho(n, k)) \ge \tau_n(\delta)$ and
$\tau_n(\rho(n, m)) \le \tau_n(r)$. Hence
$\tau_n(\delta) \le 2r + 3\delta + \tau_n(r)$, and letting
$\delta \to 0$ gives $\kappa \le 2r + \tau_n(r)$.
\end{proof}

The hypotheses hold for the toll and thin spirals, and the thin
spiral, where $\tau_n(r) = (\kappa - 2r)^+$, attains the bound for
every $r$. So a descent from $n$ to a deep point $m$ pays a toll of
at least $\kappa - 2\rho(n, m)$. If, as in
the spirals, the shift is non-contractive, (i) holds and the asymptotic distance
is $\rho$, then the toll is the defect, and the lower bound of
Theorem~\ref{thm:band} with $k = 0$ already gives
$\D(n, m) \ge c(n) - 2\rho(n, m) \ge \kappa - 2\rho(n, m)$ for deep
$m$, since $c(n) \ge \kappa$, $c(m) = 0$ and
$\lams(n, m) = 2\rho(n, m)$. The proposition obtains the same bound
from the triangle inequality alone, without the dynamics.
Remark~\ref{rem:envelopes} places the three spirals among all
defects with their descent defects.

\subsection{The ascent spiral: an ascent can carry defect}
\label{sec:ascent}

\begin{example}[Ascent spiral]\label{ex:ascent}
Let
\[
d(n, m) = \rho(n, m) +
\begin{cases}
h(2n - m) & \text{if } n < m,\\
h(n) - h(m) & \text{if } m \le n,
\end{cases}
\qquad h(n) = 2^n .
\]
\end{example}

\begin{proposition}\label{prop:ascent}
The ascent spiral is a $T_1$ quasi-metric on $\Z$ satisfying (i). The
shift is non-contractive and strictly expands every ascent. The
descent defect is $c(n) = 2^n$ and the ascent defect is $A(n) = 2^{n-1}$.
\end{proposition}

\begin{proof}
Shifting doubles both parts of the toll. Hence the shift is
non-contractive, and it strictly expands every ascent, where the toll
is positive. The triangle inequality reduces to subadditivity of the
toll, which we check case by case according to the order of $n$, $m$
and $k$. Two cases are
not immediate:
\begin{itemize}
\item an ascent $n < k$ routed through $m < n$, where
$h(2n - k) \le h(n - 1) \le h(n) - h(m)$, because $2n - k \le n - 1$
and $m \le n - 1$;
\item an ascent $n < k$ routed through $m > k$, where
$h(2n - k) \le h(m - 1) \le h(m) - h(k)$, because $2n - k \le m - 1$
and $k \le m - 1$.
\end{itemize}
The remaining cases are immediate: for an ascent $n < k$ through
$n < m < k$, use $h(2n - k) \le h(2n - m)$; for a descent $k \le n$
through $m \le k$, through $k \le m \le n$, or through $m > n$, use
$h(n) - h(k) \le h(n) - h(m)$, the identity
$h(n) - h(k) = (h(n) - h(m)) + (h(m) - h(k))$, or
$h(n) - h(k) \le h(m) - h(k)$, respectively.

For (i), pass to a subsequence with convergent base and monotone
indices. If the indices are strictly increasing, pass to a
further subsequence with $u_{i+1} \ge 2u_i + i$ for all $i$. Then
$2u_i - u_k \le -i$ whenever $i < k$, so the ascent toll
$h(2u_i - u_k)$ is at most $2^{-i}$.
If they are nonincreasing, the heights $h(u_k)$ form a monotone
bounded sequence, so the descent tolls $h(u_i) - h(u_k)$ are
eventually small.

\emergencystretch=2em
The asymptotic distance is $\rho$, since each backward step halves
both parts of the toll. Hence
$\binoppenalty=10000 \D(n, n - \Delta) = h(n) - h(n - \Delta)$ and
$\binoppenalty=10000 \D(n, n + j) = h(n - j)$, with suprema $2^n$
and $2^{n - 1}$.
\end{proof}

By Proposition~\ref{prop:afee}(a) the ascent defect is at most the
descent defect; in the ascent spiral their ratio is $1/2$ at every
point. The example refutes the guess, natural after the toll and
thin spirals, that only descents carry defect. This does not contradict
Section~\ref{sec:prelim}: there, chains make the distances
$d(b, F^\Delta b)$ of ascents, and hence their defects, small only
along the differences of a chain. The height and ascent spirals have the
same asymptotic distance and the same descent defect $c(n) = 2^n$, but different
defects on ascents, so the defect is not determined by the asymptotic distance
and the descent defect.

\begin{remark}[Mixtures]\label{rem:mix}
Over a fixed base and shift, the quasi-metrics $d = \rho + \tau$ for
which the shift is non-contractive and (i) holds form a convex set.
Indeed, a convex combination of two of them is a quasi-metric, the
shift is non-contractive for it, and every sequence has a subsequence
that is left K-Cauchy for the first, and a further subsequence that is
left K-Cauchy for the second. Mixtures of the toll and thin spirals
therefore have the same asymptotic distance and descent defect, and on each descent
from a recent to a deep point their defects fill the whole interval
between the defects of the two spirals. The set is also closed under
finite pointwise maxima, by the argument in the proof of
Proposition~\ref{prop:admissible}.
We have not formalized this remark.
\end{remark}

\needspace{12\baselineskip}
\subsection{Tolled extensions with descent zones}\label{sec:tolled}

Call a set of descents on $\Z$ a \emph{descent zone} if the shift maps
it into itself, and say that $F$ is a \emph{tolled extension with
descent zones} on an orbit if (a) of the next proposition holds. On a
single orbit, the ascents decide whether this is the case.

\begin{proposition}[Tolled extensions on an orbit]\label{prop:tolled}
Let $F$ be a non-contractive bijection of a quasi-pseudo-metric space
$(X, d)$, and identify an infinite orbit with $\Z$, so that $F$ acts on
it as the shift. The following are equivalent:
\begin{enumerate}
\item[(a)] on the orbit, $d = \ell + \sum_i \beta_i \mathbf{1}_{Z_i}$, where
$\ell$ is a shift-invariant quasi-pseudo-metric, each $\beta_i$ is
nonnegative, and the $Z_i$ are descent zones, possibly infinitely
many;
\item[(b)] the ascents of the orbit carry no defect.\enlargethispage{\baselineskip}
\end{enumerate}
In that case $\ell$ can be taken to be the asymptotic distance.
\end{proposition}

\begin{proof}
Suppose (a). Then $\ell \le d$, and since $\ell$ is shift-invariant,
$\ell(x, y) = \ell(G^k x, G^k y) \le d(G^k x, G^k y)$ for every $k$, so
$\ell \le \lam$. No zone contains an ascent, so on an ascent
$\ell \le \lam \le d = \ell$, and the defect vanishes.

Conversely, suppose (b). For each displacement $\delta \ge 1$ and each
$j \in \Z$, the set $Z_{\delta, j} = \{(n, n - \delta) : n \ge j\}$ is
a descent zone. Charge it
$\D(j, j - \delta) - \D(j - 1, j - 1 - \delta)$, which is nonnegative
by Lemma~\ref{lem:register}(b). The zones containing a descent
$(n, n - \delta)$ are the $Z_{\delta, j}$ with $j \le n$, and their
charges add up to $\D(n, n - \delta)$, because
$\D(j, j - \delta) \to 0$ as $j \to -\infty$, by the definition of
the asymptotic distance. Since ascents carry no defect,
$d = \lam + \D$ has the form (a), with $\ell = \lam$.
\end{proof}

In particular, the height and thin spirals are tolled extensions with
descent zones, and the ascent spiral is not. The substance of the
proposition is that the zones can be taken to consist of descents
only. With the sets
$Z_{\delta, j}$ for every $\delta \ne 0$, the same telescoping writes
every non-contractive bijection, on each infinite orbit, as its
asymptotic distance plus a nonnegative combination of sets of pairs that the
shift maps into themselves.

The proposition concerns the distance, not just the topology. The
identity map from the ascent spiral to the height spiral commutes
with the shift and is a homeomorphism for both topologies, although
only the height spiral is a tolled extension with descent zones.
Indeed, the two distances agree on descents, and on an ascent from
$n$ to $m > n$ they differ by $h(2n - m)$. For fixed $n$ and
$\varepsilon > 0$, this difference is at least $\varepsilon$ only for
finitely many ascents from $n$ and finitely many ascents to $n$, and
these have positive base distance. So every ball of either distance
around $n$, forward or backward, contains a ball of the other around
$n$.

\begin{proposition}[The ascent spiral is not a tolled extension]
\label{prop:bilip}
Let $(\Z, d)$ be the ascent spiral, with shift $S$. There is no
quasi-pseudo-metric space $(X', d')$ with a bijection $F'$, and no
bijection $\sigma$ from $\Z$ onto an orbit of $F'$ with
$\sigma \circ S = F' \circ \sigma$ and
$L^{-1} d(n, m) \le d'(\sigma n, \sigma m) \le L\, d(n, m)$ for some
$L \ge 1$, such that $d'$ has, on that orbit, the form (a) of
Proposition~\ref{prop:tolled}. In particular, the ascent spiral is not
conjugate to a tolled extension with descent zones by a bijection
that preserves the distance.
\end{proposition}

\begin{proof}
In the form (a), no descent zone contains an ascent, so
$d'(y, F'y) = \ell(y, F'y)$ for every $y$ on the orbit, and this
number does not depend on $y$, because $\ell$ is shift-invariant.
Since $\sigma(n + 1) = F'\sigma(n)$, the bi-Lipschitz bound gives
$d(n, n + 1) \le L\, d'(\sigma n, F'\sigma n)$, so the distances
$d(n, n + 1)$ would be bounded. In the ascent spiral, however,
$d(n, n + 1) = \rho(0, 1) + 2^{n - 1}$, which is unbounded.
\end{proof}

\subsection{Tolls over a shadow}\label{sec:tollshadow}

We now name the conclusion of Theorem~\ref{thm:structure}. The
resulting notion is a structural reading of
\cite[Question~7.6]{Duzer2026}, abstracted from the toll spiral and
the examples above.

\begin{definition}[Toll over a shadow]\label{def:toll}
A non-contractive bijection $F$ of a quasi-pseudo-metric space
$(X, d)$, with inverse $G$, is a \emph{toll over a shadow} if there is a
quasi-pseudo-metric $\ell \le d$ for which $F$ is an isometry, with
totally bounded symmetrization $\ell_s(x, y) = \ell(x, y) + \ell(y, x)$
on every orbit, such that the
defect $D_\ell = d - \ell$ and the descent defect
$c_\ell(x) = \sup_{\Delta \ge 1} D_\ell(x, G^{\Delta} x)$
satisfy, for all $x, y$ and $k \ge 0$,
\[
c_\ell(G^k x) - c_\ell(G^k y) - \ell_s(x, y) \;\le\;
D_\ell(x, y) \;\le\; c_\ell(x),
\]
with $c_\ell$ finite, $c_\ell(x) \le c_\ell(Fx)$ and
$c_\ell(G^k x) \to 0$ as $k \to \infty$ for every $x$. Any such $\ell$ is a \emph{witness}, and the \emph{shadow} is the
quasi-pseudo-metric space $(X, \ell)$, on which $F$ acts as an
isometry. $F$ is \emph{trivially tolled} if some
witness has zero defect.
\end{definition}

By Theorem~\ref{thm:structure} and Lemma~\ref{lem:register}(d), every non-contractive bijection of a space
satisfying (i) is a toll over a shadow; this is
Corollary~\ref{cor:tollanswer} below. The two-sided bound in the
definition does not determine the defect: the height and ascent spirals have the same
asymptotic distance and the same descent defect, but different defects. The
three examples above all satisfy it, in quite different ways: in the height spiral the
toll depends on the target, in the thin spiral it decreases as the
base distance grows, and in the ascent spiral ascents are charged
too.
Read as asking for a tolled extension with descent zones, up to a
conjugacy that preserves the distance or is bi-Lipschitz, the
question has a negative answer (Proposition~\ref{prop:bilip}; see
also Remark~\ref{rem:readings}). Up to a homeomorphism commuting with
the shift, the ascent spiral is not a counterexample, since the
identity is such a homeomorphism onto the height spiral (see the
discussion after Proposition~\ref{prop:tolled}).

\begin{corollary}[The structural answer to {\cite[Question~7.6]{Duzer2026}}]
\label{cor:tollanswer}
Every non-contractive bijection of a quasi-pseudo-metric space
satisfying (i), in particular of a hereditarily precompact
quasi-metric space, is a toll over the shadow $(X, \lam)$, which, on
each orbit, maps densely into the completed orbit, a compact monothetic
group on which $F$ acts as a rotation; and it is trivially tolled exactly when it is an
isometry.
\end{corollary}

\begin{proof}
By Theorem~\ref{thm:structure} and Lemma~\ref{lem:register}(d),
the asymptotic distance $\lam$ is a witness.
It is the only witness (Remark~\ref{rem:witness}), so the shadow
$(X, \lam)$ is determined by $F$, and Section~\ref{sec:shadow} shows
that on each orbit it maps densely into the completed orbit, a compact
monothetic group on which $F$ acts as a rotation. If some witness $\ell$ has
zero defect, then $\ell = d$, so $F$ is an isometry. Conversely, if $F$ is an
isometry, every ladder is constant and $\lam = d$, so the witness
$\lam$ has zero defect.
\end{proof}

\begin{remark}[The sets $U^\ell_\varepsilon$, isometries without a witness, and uniqueness of the witness]
\label{rem:witness}
Defects of at least $\varepsilon$ occur only on pairs whose source
lies in the set
$U^{\ell}_\varepsilon = \{x : c_\ell(x) \ge \varepsilon\}$. Each such set satisfies $F(U^{\ell}_\varepsilon) \subseteq
U^{\ell}_\varepsilon$, since $c_\ell$ is nondecreasing along orbits,
and every backward orbit eventually leaves it, since $c_\ell$
tends to $0$ along backward orbits. For the toll spiral and $0 < \varepsilon
\le \kappa$, the set $U_\varepsilon$ is the set of recent integers, and
the toll zone of \cite[Remark~5.8]{Duzer2026} is the set of descents
whose source lies in it. Without (i), an isometry need not be a
toll over a shadow at all: the shift on $(\Z, |n - m|)$ is an
isometry with no witness. If $\ell$ were a witness, its
symmetrization would be totally bounded on the orbit, so $\ell$
would be bounded there. Then $D_\ell(n, n - \Delta) \ge \Delta -
\sup \ell$ would be unbounded in $\Delta$, and $c_\ell$ would be
infinite, whereas a witness has $c_\ell$ finite. Under (i) this cannot happen: by
Theorem~\ref{thm:structure} and Lemma~\ref{lem:register}(d), the
asymptotic distance $\lam$ is always a witness. It is also the only one. Let
$\ell$ be a witness. Since $F$ is an isometry for $\ell$ and
$\ell \le d$, we have $\ell(x, y) = \ell(G^k x, G^k y) \le
d(G^k x, G^k y)$ for every $k$, so $\ell \le \lam$. Conversely,
$d(G^k x, G^k y) - \ell(x, y) = D_\ell(G^k x, G^k y) \le
c_\ell(G^k x)$, which tends to $0$, so $\lam \le \ell$. For
$\ell = \lam$ we have $D_\ell = \D$, $c_\ell = c$ and
$\ell_s = \lams$, so for this witness the notation of
Definition~\ref{def:toll} reduces to that of
Definition~\ref{def:fee}. Finally, the inequality $\ell \le \lam$
used only that $\ell \le d$ and that $F$ is an isometry for $\ell$,
and not (i) or the two-sided bound. Since $\lam$ itself has these two properties
(Lemma~\ref{lem:register}(a)), it is the largest quasi-pseudo-metric
below $d$ for which $F$ is an isometry.
\end{remark}

\section{The descent defect tends to zero along backward orbits}\label{sec:depth}

By Proposition~\ref{prop:fee}(c), the descent defect is nondecreasing along
orbits, so its backward limit $c_\infty(a) = \lim_k c(G^k a)$ exists.
We show that this limit is always $0$. Only hypothesis (i) is used.

The proof rests on Lemma~\ref{lem:feetarget}. It says that every
point $x$ has a backward iterate $y = G^\Delta x$ for which
$\D(x, y)$ is close to $c(x)$, while $\lams(x, y)$ and $\D(y, x)$
are small; it follows from Theorem~\ref{thm:source} and two more
uses of (i). To prove Theorem~\ref{thm:vanish}, we suppose that the
descent defect stays above some $\eta > 0$ along a backward orbit.
Repeated use of the lemma then gives a sequence $(u_m)$ with
$d(u_i, u_j) > \eta/2$ whenever $i < j$, and no subsequence of it is
left K-Cauchy, contrary to (i).

The theorem says nothing about how fast the descent defect tends
to $0$. In the height spiral, the descent defect
$c(n) = 2^n$ never reaches $0$, while for zone tolls it is $0$ from
some depth on. Both cases are treated after
the vanishing theorem (Remark~\ref{rem:rate} and
Theorem~\ref{thm:zone}), and the section ends with the proof of
Theorem~\ref{thm:structure}.

\subsection{The vanishing theorem}

\begin{lemma}[The descent defect is nearly attained on near-returns]\label{lem:feetarget}
Let $x \in X$, $\eta \le c(x)$ and $\delta > 0$. There is $\Delta
\ge 1$ with
\[
\D(x, G^\Delta x) > \eta - \delta, \qquad
\lams(x, G^\Delta x) < \delta, \qquad
\D(G^\Delta x, x) < \delta .
\]
\end{lemma}

\begin{proof}
Since $\eta \le c(x)$, there is $b = G^{\Delta_0} x$ with
$\D(x, b) > \eta - \delta/4$. Since the ladder of $(x, b)$ decreases
to $\lam(x, b)$, there is $K$ with
$d(G^k x, G^k b) < \lam(x, b) + \delta/4$ for $k \ge K$.
Theorem~\ref{thm:source}, applied with $\eta' = \D(x, b) - \delta/4$,
$M = K$ and tolerance $\delta/4$, gives a chain $(t_i)$ along whose
differences $\Delta$
\[
d(x, G^\Delta x) > \D(x, b) - \delta/2 .
\]
Apply (i) to the sequence $(G^{t_i} x)_i$ and then, along the
subchain obtained, to $(F^{t_i} x)_i$. After discarding finitely many
terms, this gives a subchain, still denoted $(t_i)$, with
$d(G^{t_i} x, G^{t_j} x) < \delta/4$ and
$d(F^{t_i} x, F^{t_j} x) < \delta/4$ for $i < j$.

Fix $i < j$ and put $p = t_i$, $q = t_j$ and $\Delta = q - p$. By
invariance, $\lam(x, G^\Delta x) \le d(G^p x, G^q x) < \delta/4$.
Since $F^q$ maps $(G^\Delta x, x)$ to $(F^p x, F^q x)$, invariance
also gives
\[
\lam(G^\Delta x, x) = \lam(F^p x, F^q x) \le d(F^p x, F^q x) <
\delta/4 ,
\]
and non-contractivity of $F^q$ gives
\[
\D(G^\Delta x, x) \le d(G^\Delta x, x) \le d(F^p x, F^q x) <
\delta/4 .
\]
The first two bounds give
$\lams(x, G^\Delta x) < \delta/2$, and the first, together with the
display, gives
$\D(x, G^\Delta x) > \D(x, b) - 3\delta/4 > \eta - \delta$.
\end{proof}

\begin{theorem}[Decay of the descent defect]\label{thm:vanish}
For every $a \in X$, $c(G^k a) \to 0$ as $k \to \infty$.
\end{theorem}

\begin{proof}
By Proposition~\ref{prop:fee}(c), the descent defects $c(G^k a)$ are
nonincreasing in $k$. Suppose, for a contradiction, that
$c(G^k a) \ge \eta > 0$ for all $k$. Put $\varepsilon = \eta/16$
and $w_m = \varepsilon 2^{-m}$.

We build a sequence $(u_m)$ on the backward orbit of $a$. Let
$u_0 = a$ and, given $u_m$, let $u_{m+1} = G^{\Delta_m} u_m$, where
$\Delta_m$ is given by Lemma~\ref{lem:feetarget} for $x = u_m$ and
$\delta = w_m$; this is possible because $c(u_m) \ge \eta$. Thus
\[
\D(u_m, u_{m+1}) > \eta - w_m, \qquad
\lams(u_m, u_{m+1}) < w_m, \qquad
\D(u_{m+1}, u_m) < w_m .
\]
Since $\sum_{m=i}^{j-1} w_m = 2w_i - 2w_j$, the triangle inequality
for $\lams$ gives $\lams(u_i, u_j) \le 2w_i - 2w_j$ for $i \le j$.
Likewise, Lemma~\ref{lem:register}(c), applied one step at a time
with the last two bounds, gives $\D(u_j, u_i) \le 4w_i - 4w_j$ for
$i \le j$.

Now let $i < j$. Lemma~\ref{lem:register}(c), applied to the points
$u_{j-1}$, $u_i$, $u_j$, gives
\[
\D(u_i, u_j) \ge \D(u_{j-1}, u_j) - \D(u_{j-1}, u_i)
- \lams(u_{j-1}, u_i)
> \eta - w_{j-1} - 4\varepsilon - 2\varepsilon
\ge \eta - 7\varepsilon > \eta/2 .
\]
Hence $d(u_i, u_j) \ge \D(u_i, u_j) > \eta/2$ for all $i < j$, so no
subsequence of $(u_m)$ is left K-Cauchy, contradicting (i).
\end{proof}

\begin{corollary}[Asymptotic isometry along backward orbits]\label{cor:rigid}
For every $a \in X$ and $\varepsilon > 0$ there is $k$ such that
$d(G^{k'} a, y) < \lam(G^{k'} a, y) + \varepsilon$ for all $k' \ge
k$ and all $y \in X$. In particular,
$0 \le d(Fx, Fy) - d(x, y) < \varepsilon$ for $x = G^{k'} a$ with
$k' > k$ and every $y$: the map is asymptotically an isometry along
its backward orbits. The ascent defect also tends to $0$ along
backward orbits.
\end{corollary}

\begin{proof}
By Proposition~\ref{prop:fee}(c), the descent defects $c(G^k a)$ are
nonincreasing in $k$, so Theorem~\ref{thm:vanish} gives $k$ with
$c(G^{k'} a) < \varepsilon$ for all $k' \ge k$. By
Proposition~\ref{prop:fee}(b), $\D(G^{k'} a, y) \le c(G^{k'} a) <
\varepsilon$ for every $y$, which is the first claim. For the
second, let $x = G^{k'} a$ with $k' > k$. Then $Fx = G^{k'-1} a$ with
$k' - 1 \ge k$, so the first claim, invariance and $\lam \le d$ give
\[
d(Fx, Fy) < \lam(Fx, Fy) + \varepsilon = \lam(x, y) + \varepsilon
\le d(x, y) + \varepsilon;
\]
the lower bound is non-contractivity. Finally, the ascent defect
tends to $0$ because $A \le c$ (Proposition~\ref{prop:afee}(a)).
\end{proof}

It follows that every backward orbit is totally bounded for the
symmetrized distance $d(x, y) + d(y, x)$. Given $\varepsilon > 0$,
take $k$ as in the corollary. For $k', k'' \ge k$, the symmetrized
distance between $G^{k'} a$ and $G^{k''} a$ is less than
$\lams(G^{k'} a, G^{k''} a) + 2\varepsilon$, and $\lams$ is totally
bounded on the orbit by Theorem~\ref{thm:structure}(a). So the tail
$\{G^{k'} a : k' \ge k\}$ contains finitely many points such that
every point of the tail is within $2\varepsilon$ of one of them for
$\lams$, and hence within $4\varepsilon$ for the symmetrized distance.
Adding the points $G^{k'} a$ with $k' < k$ gives a finite
$4\varepsilon$-net for the whole backward orbit.

Moreover, there are arbitrarily large $n$ for which $G^n$ moves every
point of the tail by less than $2\varepsilon$ in both directions. Let $k$
be as above and $N \ge 1$. Theorem~\ref{thm:shadow}, applied to the
sequence $n_i = -iN$, gives $i < j$ with
$\lams(G^{iN} a, G^{jN} a) < \varepsilon$. By invariance,
$n = (j - i)N$ is a multiple of $N$ with $n \ge N$ and
$\lams(a, G^n a) < \varepsilon$. Let $x = G^{k'} a$ with $k' \ge k$.
Then $G^n x$ also lies in the tail, so the corollary and the
invariance of $\lam$ give
\[
d(x, G^n x) < \lam(a, G^n a) + \varepsilon < 2\varepsilon,
\qquad
d(G^n x, x) < \lam(G^n a, a) + \varepsilon < 2\varepsilon .
\]
Neither property
need hold for forward orbits. In the toll spiral, the forward orbit
of $0$ consists of the recent integers, and every pair of distinct
recent integers is at symmetrized distance at least $\kappa$, since
one of the two pairs is a descent from a recent source. So this
orbit is not totally bounded, and $d(F^n x, x) \ge \kappa$ for every
recent $x$ and every $n \ge 1$.

With the forward returns of \cite[Lemma~3.3]{DuzerC}, these
near-returns show that a non-isometric example cannot be
placed inside a space with limits. For a sequence $(x_n)$ in a
quasi-pseudo-metric space $(Z, e)$, a point $y$ is an \emph{upper
limit} if $e(x_n, y) \to 0$, and a point $p$ is a \emph{lower limit}
if $e(p, x_n) \to 0$; Yoneda limits and sequential Yoneda
completeness are as in \cite[Definition~2.1]{DuzerC}, and $R_0$ is
as in Section~\ref{sec:prelim}.

\begin{corollary}[No $R_0$ host with limits]\label{cor:host}
Suppose that $F$ is not an isometry, and that $(X, d)$ is a subspace
of an $R_0$ quasi-pseudo-metric space $(Z, e)$. Then some left
K-Cauchy sequence of $X$ has no upper limit in $Z$, and some left
K-Cauchy sequence of $X$ has no lower limit in $Z$. In particular, no
sequentially Yoneda complete quasi-metric space, and no quasi-metric
space whose topology or conjugate topology is countably compact,
contains $(X, d)$.
\end{corollary}

\begin{proof}
\emph{Upper limits.} We follow the proof of
\cite[Theorem~3.4]{DuzerC}. Suppose that every left K-Cauchy sequence
of $X$ has an upper limit in $Z$, and let $a, b \in X$. By
\cite[Lemma~3.3]{DuzerC}, which uses only (i), there are $m_i \ge 1$
with $d(a, F^{m_i} a) \to 0$ and $d(b, F^{m_i} b) \to 0$. By (i),
after passing to a subsequence, $(F^{m_i} a)$ is left K-Cauchy; let $y
\in Z$ be an upper limit of it. Then
$e(a, y) \le d(a, F^{m_i} a) + e(F^{m_i} a, y) \to 0$, so $e(a, y) =
0$, hence $e(y, a) = 0$, and $d(F^{m_i} a, a) \le e(F^{m_i} a, y) +
e(y, a) \to 0$. Since $F^{m_i - 1}$ is non-contractive,
\[
d(a, b) \le d(Fa, Fb) \le d(F^{m_i} a, F^{m_i} b) \le
d(F^{m_i} a, a) + d(a, b) + d(b, F^{m_i} b),
\]
and the right side tends to $d(a, b)$. So $F$ would be an isometry.

\emph{Lower limits.} Since $F$ is not an isometry, some pair has
$d(Fx, Fy) > d(x, y)$, hence $\D(Fx, Fy) > \D(x, y) \ge 0$ by the invariance of $\lam$
(Lemma~\ref{lem:register}(a)), and $a = Fx$ has $c(a) > 0$ by
Proposition~\ref{prop:fee}(b). By Theorem~\ref{thm:shadow}, as in
the second paragraph after Corollary~\ref{cor:rigid}, there are
$t_i \to \infty$ with $\lams(a, G^{t_i} a) \to 0$, and by (i)
we may assume that $(G^{t_i} a)$ is left K-Cauchy. Suppose that it has
a lower limit $p \in Z$. By Proposition~\ref{prop:fee}(b) and
Theorem~\ref{thm:vanish},
\[
d(G^{t_i} a, a) \le \lam(G^{t_i} a, a) + c(G^{t_i} a) \to 0,
\]
so $e(p, a) \le e(p, G^{t_i} a) + d(G^{t_i} a, a) \to 0$. Hence
$e(p, a) = 0$, so $e(a, p) = 0$, and $d(a, G^{t_i} a) \le e(a, p) +
e(p, G^{t_i} a) \to 0$. But by Theorem~\ref{thm:band} with $k = 0$,
\[
d(a, G^{t_i} a) \ge \D(a, G^{t_i} a) \ge c(a) - c(G^{t_i} a) -
\lams(a, G^{t_i} a) \to c(a) > 0 .
\]

Finally, Yoneda limits are upper limits (take $z = y$ in their
definition). If the topology of $e$ is countably compact, every left
K-Cauchy sequence has a lower limit, and if the conjugate topology is,
every left K-Cauchy sequence has an upper limit
\cite[Lemma~2.3(a)--(c)]{DuzerC}, with (a) applied to $e$ and to
$e^{-1}$.
\end{proof}

This extends \cite[Proposition~5.2]{DuzerC} from the toll spiral to
every non-isometric example. The first half uses only (i) and the
returns of \cite[Lemma~3.3]{DuzerC}, so it holds for every
non-contractive self-map of a space satisfying (i), bijective or not.
The second half uses the structure theory: the backward near-returns,
the decay of the descent defect, and the lower bound of
Theorem~\ref{thm:band}.

\begin{remark}[The rate is not forced]\label{rem:rate}
In Theorem~\ref{thm:vanish}, the descent defect tends to $0$ but need not
reach it. In the height and ascent spirals, the descent defect $c(n) = 2^n$ is
positive at every integer, and so is the descent defect $g(n)$ of the integer
layer of the space of \cite[Proposition~6.1]{DuzerC}. In fact any
rate of decay occurs: every nonconstant nondecreasing function
$g \ge 0$ tending to $0$ in the deep past is a descent defect. Charge each
descent from $n$ the amount $g(n)$ on top of the base distance. As
for the
integer layer, the result is a $T_1$ quasi-metric on $\Z$ satisfying
(i), on which the shift is non-contractive and not an isometry; its
asymptotic distance is $\rho$, and its descent defect is $g$.
Theorem~\ref{thm:pairs}, with $A = 0$, gives the same toll over every
shift-invariant quasi-pseudo-metric base satisfying (i), again with
descent defect $g$; the result is a quasi-metric when the base is one,
and the shift is not an isometry since $g$ is not constant.
Theorem~\ref{thm:pairs} also covers pairs of descent and ascent
defects, and Theorem~\ref{thm:occur} extends the realization to every
asymptotic distance. For $g = \kappa\,\mathbf{1}[n \ge 0]$ this is the
toll spiral. Zone tolls, treated next, behave differently from a
general $g$: their descent defect is always $0$ below some depth, as
for the toll spiral.
\end{remark}

\subsection{Zone tolls: the descent defect is eventually zero along backward orbits}

Zone tolls are two-valued tolls on $\Z$: a descent costs either
nothing or a fixed amount $\kappa$, according to its source and its
displacement. Formally, a \emph{displacement zone} is a set
$Z \subseteq \Z \times \N_{\ge 1}$ of pairs (source, displacement),
not a set of points like the sets $U_\varepsilon$ of
Theorem~\ref{thm:structure}, and
$Z_n = \{\delta : (n, \delta) \in Z\}$ is the set of displacements
charged from $n$. Consider the conditions
\begin{itemize}
\item[(Z1)] $(n - 1, \delta) \in Z$ implies $(n, \delta) \in Z$;
\item[(Z2)] $(n, \delta + \delta') \in Z$ implies $(n, \delta) \in Z$
or $(n - \delta, \delta') \in Z$.
\end{itemize}
Condition (Z1) says that the toll $\kappa\,\mathbf{1}_Z$ is
shift-monotone, and (Z2) that it is subadditive along the route
$n \to n - \delta \to n - \delta - \delta'$. The toll spiral is the
zone toll given by the displacement zone with $Z_n = \N_{\ge 1}$ for $n \ge 0$ and $Z_n = \emptyset$ for
$n < 0$. The next theorem adds the toll to a shift-invariant base $\rho$. By
(Z1), the shift is non-contractive for the resulting distance $d$,
so the asymptotic distance and the descent defect are defined as in
Section~\ref{sec:prelim}. The charged descents
$\{(n, n - \delta) : \delta \in Z_n\}$ form a descent zone in the
sense of Section~\ref{sec:tolled}, so when $\rho$ and $d$ are
quasi-pseudo-metrics, $d$ is a tolled extension with a single
descent zone. Conditions (Z1) and (Z2) alone do not make $d$ a
quasi-pseudo-metric, even when $\rho$ is a metric
(Remark~\ref{rem:z1z2}).

\begin{theorem}[Zone tolls]\label{thm:zone}
Let $Z$ satisfy (Z1) and (Z2), let $\rho$ be a shift-invariant
function on $\Z \times \Z$, let $\kappa > 0$, and let
\[
d(n, m) = \rho(n, m) + \kappa\,\mathbf{1}[m < n \text{ and } n - m \in Z_n].
\]
Then there is $n_0$ such that $Z_n = Z_{n_0}$ for all $n \le n_0$, and
for every $n \le n_0$ and every $m$ the ladder of $(n, m)$ is
constant; in particular every defect from a point $n \le n_0$
vanishes, so the descent defect is $0$ at every $n \le n_0$.
\end{theorem}

\begin{proof}
Let $C_n = \{0\} \cup (\N_{\ge 1} \setminus Z_n)$ be the set of
displacements that are free from $n$, together with $0$. We first
show that $C_n$ is a submonoid of $\N$. Let $\delta, \delta' \ge 1$
be free from $n$, and suppose that $\delta + \delta'$ were charged
from $n$. By (Z2), either $\delta$ would be charged from $n$, which
it is not, or $\delta'$ would be charged from $n - \delta$; in the
latter case (Z1), applied $\delta$ times, would charge $\delta'$ from
$n$, which it is not either.

By (Z1), $C_n \subseteq C_{n-1}$: the free monoids grow with depth.
An increasing chain of submonoids of $\N$ is eventually constant. Its
union is a submonoid of $\N$, so it has finitely many generators, as every submonoid of $\N$ does.
They all lie in one member of the chain, and from that member on,
the chain equals the union. Applied to the chain $k \mapsto C_{-k}$, this gives
$n_0$ with $C_n = C_{n_0}$, that is $Z_n = Z_{n_0}$, for all
$n \le n_0$.

Finally, let $n \le n_0$ and let $m$ be arbitrary. For every
$k \ge 0$, the pair $(n - k, m - k)$ has the same order and the same
displacement as $(n, m)$, and $n - k \le n_0$, so it pays the same
toll. Since $\rho$ is invariant, the ladder of $(n, m)$ is constant;
hence the asymptotic distance equals $d$, and the defect is $0$.
\end{proof}

So, for a zone toll, the charges eventually stop changing as one goes
into the past, and from then on they are part of the asymptotic distance. The
reason is the ascending chain condition for submonoids of $\N$:
charges that kept changing forever would give an infinite strictly
increasing chain of free monoids.

\begin{remark}[(Z1)--(Z2) and the triangle inequality]\label{rem:z1z2}
The proof uses (Z1), (Z2) and the invariance of $\rho$, but never
the triangle inequality for $d$, and the pair (Z1)--(Z2) is
independent of it. Condition (Z2)
is the triangle inequality for the toll, restricted to descending
routes $n \to n - \delta \to n - \delta - \delta'$. For instance, the
displacement zone $Z_n = \{1\}$ satisfies (Z1) and (Z2), but $d$ fails the triangle
inequality along the route $n \to n + 1 \to n - 1$ as soon as $\kappa$
exceeds the extra base cost of that detour,
$\rho(n, n + 1) + \rho(n + 1, n - 1) - \rho(n, n - 1)$. In the other
direction, $d$ can be a quasi-metric while (Z2) fails. Take the displacement zone
$Z_n = \{2\}$ over the
base $\rho(n, m) = \sqrt{|n - m|}$, with $0 < \kappa \le 2 - \sqrt{2}$.
Then $d$ is a quasi-metric. Only descents of displacement $2$ are
charged, and any detour around such a descent costs at least $2$ in
the base, which is at least $\sqrt{2} + \kappa$. Yet (Z2) fails for
$\delta = \delta' = 1$.
\end{remark}

\subsection{Proof of the structure theorem}

By definition of the defect, $d = \lam + \D$ for every
non-contractive bijection. Moreover, $\D = 0$ exactly when $F$ is an
isometry: if $F$ is an isometry, every ladder is constant, so
$\lam = d$; and if $\D = 0$, then $d = \lam$, for which $F$ is an
isometry. So the non-isometric case is exactly the case $\D \ne 0$.

\begin{proof}[Proof of Theorem~\ref{thm:structure}]
(a) is Lemma~\ref{lem:register}(a) and Theorem~\ref{thm:shadow};
(b) is Theorem~\ref{thm:band}, Proposition~\ref{prop:fee}(b),(c)
and Theorem~\ref{thm:vanish}; the consequences follow from (b).
\end{proof}

The examples of Section~\ref{sec:spirals} show that the bounds of
Theorem~\ref{thm:structure}(b) do not determine the defect, and
Remark~\ref{rem:rate} shows that no rate of decay is forced.

\section{Which defects occur}\label{sec:occur}

Theorem~\ref{thm:structure} constrains the defect of a map through
its asymptotic distance and its descent defect. This section asks
which defects occur once the asymptotic distance is fixed. Fix a set
$X$ with a bijection $F$, with inverse $G$, and a quasi-pseudo-metric
$\ell$ on $X$ that satisfies (i) and for which $F$ is an isometry. A
function $\tau \ge 0$ on $X \times X$ is \emph{admissible over
$\ell$} if $d = \ell + \tau$ is a quasi-pseudo-metric satisfying
(i), $F$ is non-contractive for $d$, and $\ell$ is the asymptotic
distance of $d$. Then $\tau$ is
the defect of $d$, and every non-contractive bijection of a space
satisfying (i) arises in this way, with $\ell$ its asymptotic distance
and $\tau$ its defect. For any $\tau \ge 0$ we write $c_\tau(x) =
\sup_{\Delta \ge 1} \tau(x, G^\Delta x)$; for an admissible toll this
is the descent defect. The \emph{slack} of a route $x \to y \to z$
is $\ell(x, y) + \ell(y, z) - \ell(x, z) \ge 0$.

\begin{proposition}[Admissible tolls]\label{prop:admissible}
A function $\tau \ge 0$ on $X \times X$ is admissible over $\ell$ if
and only if
\begin{enumerate}
\item[(a)] $\tau(x, x) = 0$;
\item[(b)] $\tau(Fx, Fy) \ge \tau(x, y)$;
\item[(c)] $\tau(G^k x, G^k y) \to 0$ as $k \to \infty$;
\item[(d)] $\tau(x, z) \le \tau(x, y) + \tau(y, z) + \ell(x, y) +
\ell(y, z) - \ell(x, z)$, that is, $\tau$ satisfies the triangle
inequality up to the slack;
\item[(e)] $\ell + \tau$ satisfies (i).
\end{enumerate}
The admissible tolls over $\ell$ form a convex set that is closed
under finite pointwise maxima, and $\tau + \sigma$ is admissible whenever
$\tau$ is admissible and $\sigma$ is an admissible toll satisfying the
triangle inequality.
\end{proposition}

\begin{proof}
Conditions (a) and (d) say that $\ell + \tau$ is a quasi-pseudo-metric,
and (b) that $F$ is non-contractive for it, since $\ell$ is
$F$-invariant. Given (b), the ladder of $(x, y)$ is
$\ell(x, y) + \tau(G^k x, G^k y)$, nonincreasing in $k$, and its
limit is $\ell(x, y)$ exactly when (c) holds. Conditions (a)--(d)
pass to convex combinations, since (d) is affine in $\tau$ with a
slack that does not depend on $\tau$; to maxima, since the
inequalities (d) for $\tau_1$ and $\tau_2$ give it for
$\max(\tau_1, \tau_2)$; and to $\tau + \sigma$, adding (d) for $\tau$
to the triangle inequality for $\sigma$. For (e), each new distance is
at most $\max(\ell + \tau_1, \ell + \tau_2)$ or $(\ell + \tau) + (\ell
+ \sigma)$, and a subsequence that is left K-Cauchy for the first
distance has a further subsequence that is left K-Cauchy for the
second.
\end{proof}

Condition (e) has a combinatorial form.

\begin{theorem}[Condition (i) for $\ell + \tau$]\label{thm:ramsey}
Let $\tau \ge 0$. Then $\ell + \tau$ satisfies (i) if and only if
no sequence $(x_i)$ has $\inf_{i < j} \tau(x_i, x_j) > 0$.
\end{theorem}

\begin{proof}
A sequence with $\tau(x_i, x_j) \ge \varepsilon$ for all $i < j$ has no
left K-Cauchy subsequence for $\ell + \tau$, since its subsequences
have the same property. Conversely, applied to constant sequences,
the hypothesis gives $\tau(x, x) = 0$. Let $(x_i)$ be a sequence; by (i)
for $\ell$ we may assume that it is left K-Cauchy for $\ell$. Put $I_0
= \N$. For $m \ge 1$, color each pair $i < j$ in $I_{m-1}$ according
to whether $\tau(x_i, x_j) \ge 1/m$. By the infinite Ramsey theorem
there is an infinite $I_m \subseteq I_{m-1}$ on which all pairs have
the same color, and by hypothesis the color is ``$< 1/m$''. Choose
$i_1 < i_2 < \cdots$ with $i_m \in I_m$. For $m < m'$ both indices lie
in $I_m$, so $\tau(x_{i_m}, x_{i_{m'}}) < 1/m$. Since $(x_{i_m})$ is
also left K-Cauchy for $\ell$, it is left K-Cauchy for $\ell + \tau$.
\end{proof}

\subsection{Every descent defect occurs}

For $c : X \to [0, \infty)$ let $\tau_c$ be as in
Section~\ref{sec:band}. The next theorem is the converse of the
statements on $c$ in Theorem~\ref{thm:structure}(b). Together with
Proposition~\ref{prop:realize}, it shows that neither the completed
orbits nor the descent defect are constrained beyond what
Theorems~\ref{thm:structure} and~\ref{thm:shadtop} state.

\begin{theorem}[Every descent defect occurs]\label{thm:occur}
A function $c : X \to [0, \infty)$ is the descent defect of an
admissible toll over $\ell$ if and only if $c(x) \le c(Fx)$ and
$c(G^k x) \to 0$ as $k \to \infty$, for every $x$. In that case
$\tau_c$ is admissible, satisfies the triangle inequality, and has
descent defect $c$; and $\tau_c \le \tau$ for every admissible toll
$\tau$ with descent defect $c$ that satisfies the triangle inequality.
\end{theorem}

The proof uses a combinatorial fact. Order the nonincreasing
sequences $N : \N \to \N$ with finite support pointwise; these are
the Young diagrams, ordered by inclusion. This order is a
well-quasi-order: every sequence $(N_i)$ has $i < j$ with $N_i \le
N_j$. Indeed, write $N$ as the finite word of its nonzero values
$N(0) \ge N(1) \ge \dots \ge N(r)$. By Higman's lemma, finite words
over the well-quasi-order $\N$ are well-quasi-ordered by the
subword order. If the word of $N$ embeds into that of $N'$, with
$N(q) \le N'(s_q)$ and $s_0 < s_1 < \cdots$, then $s_q \ge q$, and
$N(q) \le N'(s_q) \le N'(q)$ since $N'$ is nonincreasing.

\begin{proof}
Necessity is Proposition~\ref{prop:fee}(c) and
Theorem~\ref{thm:vanish}, applied to $(X, \ell + \tau)$. Conversely,
let $c$ satisfy the two conditions. Each function
$(c(G^k x) - c(G^k y))^+$ satisfies the triangle inequality, so
$\tau_c$ does, and $\tau_c(x, x) = 0$; this gives (a) and (d).
Since $G^k F = G^{k-1}$, the supremum defining $\tau_c(Fx, Fy)$ runs
over $k \ge -1$, which gives (b). Since $c(G^k x) \le c(x)$ for $k \ge
0$, we have $\tau_c(x, y) \le c(x)$, and (c) follows from
$c(G^k x) \to 0$. The descent defect of $\tau_c$ is at most $c(x)$,
and at least $\tau_c(x, G^\Delta x) \ge c(x) - c(G^\Delta x)$ for every
$\Delta$, which tends to $c(x)$. For (e), let $\varepsilon > 0$ and
let $(x_i)$ be a sequence. The functions $N_i(k) = \lfloor 2 c(G^k
x_i)/\varepsilon \rfloor$ are nonincreasing in $k$, since $c$ is
nondecreasing along orbits, and have finite support, since $c(G^k x_i)
\to 0$. So there are $i < j$ with $N_i \le N_j$, and then $c(G^k x_i)
- c(G^k x_j) < \varepsilon/2$ for every $k$, so $\tau_c(x_i, x_j) \le
\varepsilon/2 < \varepsilon$. By Theorem~\ref{thm:ramsey}, (e) holds.

Finally, let $\tau$ be admissible with descent defect $c$ and satisfy
the triangle inequality. For $\Delta \ge 1$, $\tau(x, y) \ge \tau(x,
G^\Delta x) - \tau(y, G^\Delta x) \ge \tau(x, G^\Delta x) - c(y)$, by
Proposition~\ref{prop:fee}(b) applied to $(X, \ell + \tau)$. Taking
the supremum over $\Delta$ gives $\tau(x, y) \ge c(x) - c(y)$, and
applying this to $(G^k x, G^k y)$, with (b), gives $\tau \ge \tau_c$.
\end{proof}

On an orbit, the toll of Remark~\ref{rem:rate}, which charges $c(x)$
on every descent from $x$ and nothing on ascents, is admissible by
Theorem~\ref{thm:pairs} below, with $A = 0$, and attains the upper
bound on every descent. Across orbits, the analogous toll, charging
$c(x)$ on the backward iterates of $x$ and nothing elsewhere, can
violate (d): on $\Z \times \{0, 1\}$ with $\ell((n, i), (m, j)) =
\rho(n, m)$, $c = \kappa\,\mathbf{1}[n \ge 0]$ on the first copy and $c = 0$ on
the second, the route $(n, 0) \to (n, 1) \to (n - 1, 0)$ with $n \ge
0$ has slack $0$ and costs $0 < \kappa$. The toll $\tau_c$ avoids this by comparing
two backward orbits step by step.

\subsection{On one orbit}

Now let $X = \Z$, let $F$ be the shift, and let $\ell = \rho$ be a
shift-invariant quasi-pseudo-metric satisfying (i), for instance the
arc metric of Section~\ref{sec:spirals}. Write $\rho_s(\delta) =
\rho(n, n - \delta) + \rho(n - \delta, n)$, which does not depend on
$n$. On one orbit, condition (e) splits into a condition on descents
and a condition on ascents.

\begin{corollary}[Condition (i) on an orbit]\label{cor:chains}
Let $\tau \ge 0$ satisfy (a)--(d) of
Proposition~\ref{prop:admissible}. Then $\rho + \tau$ satisfies (i)
if and only if
\begin{enumerate}
\item[(V)] $c_\tau(n) \to 0$ as $n \to -\infty$, and
\item[(S)] no sequence $n_1 < n_2 < \cdots$ has $\inf_{i < j}
\tau(n_i, n_j) > 0$.
\end{enumerate}
\end{corollary}

\begin{proof}
If (i) holds, then $(\Z, \rho + \tau)$ is a space as in
Section~\ref{sec:prelim}, with asymptotic distance $\rho$ by (c) and
defect $\tau$; so (V) is Theorem~\ref{thm:vanish}, and (S) holds by
Theorem~\ref{thm:ramsey}. Conversely, suppose that $\tau(n_i, n_j)
\ge \varepsilon$ for all $i < j$. Passing to a subsequence, $(n_i)$ is
constant or strictly monotone. It is not constant, since $\tau(n, n) =
0$. If it is strictly decreasing, then $n_i \to -\infty$ and $\tau(n_i,
n_j) \le c_\tau(n_i) \to 0$ by (V). So it is strictly increasing, which
(S) excludes. Theorem~\ref{thm:ramsey} concludes.
\end{proof}

Condition (S) concerns ascents only. For tolls that vanish on ascents
it holds trivially, and Corollary~\ref{cor:chains} reduces (i) to (V).

\begin{theorem}[Admissible tolls on an orbit]\label{thm:orbit}
Let $c : \Z \to [0, \infty)$ be nondecreasing with $c(n) \to 0$ as
$n \to -\infty$. A function $\tau \ge 0$ on $\Z \times \Z$ is an
admissible toll over $\rho$ with descent defect $c$ if and only if it
satisfies (a), (b) and (d) of Proposition~\ref{prop:admissible},
condition (S), and
\[
\bigl(\tau_c(n, n - \delta) - \rho_s(\delta)\bigr)^+ \le \tau(n, n -
\delta) \le c(n) \qquad (n \in \Z,\ \delta \ge 1).
\]
\end{theorem}

\begin{proof}
Necessity follows from Proposition~\ref{prop:admissible},
Corollary~\ref{cor:chains} and the bounds of Section~\ref{sec:band}.
For sufficiency we first derive (c). By (b), the limit $E(x, y) =
\lim_{k \to \infty} \tau(x - k, y - k)$ exists; $E \le \tau$, $E$ is
shift-invariant, and it satisfies (d). On descents $E = 0$, since
$\tau(n, n - \delta) \le c(n)$ and $c \to 0$ in the past. Suppose that
$E(x, x + q) = L > 0$ for some $q \ge 1$, and let $0 < \eta < L$. For
$r > q$, the route $x \to y = x + q - r \to x + q$ starts with a
descent, and its slack $\rho(x, y) + \rho(y, x + q) - \rho(x, x + q)$
is at most $\rho_s(r)$, since $\rho(x, y) \le \rho(x, x + q) + \rho(x
+ q, y)$. So $E(n, n + r) \ge L - \rho_s(r)$ for all $n$. Color the
pairs $i < j$ of natural numbers according to whether $\rho_s(j - i) <
\eta$. By Theorem~\ref{thm:shadow}, applied to $(\Z, \rho)$, for
which the shift is an isometry and the asymptotic distance is $\rho$, no
infinite set has all its pairs of the second color, so by Ramsey's
theorem some infinite set has all its pairs of the first color; thin
it to $n_1 < n_2 < \cdots$ with gaps larger than $q$. Then $\tau(n_i,
n_j) \ge E(n_i, n_j) \ge L - \eta$ for all $i < j$, contradicting (S).
So $E = 0$, which is (c).

The upper bound gives (V), so $\tau$ is admissible by
Corollary~\ref{cor:chains}. Its descent defect is at most $c$ by the
upper bound. It is at least $c$: given $n$ and $\varepsilon > 0$,
choose $M$ with $c(n - M) < \varepsilon$, and, by
Theorem~\ref{thm:shadow} as in the second paragraph after
Corollary~\ref{cor:rigid}, some $\delta \ge M$ with $\rho_s(\delta) <
\varepsilon$; then $\tau(n, n - \delta) \ge c(n) - c(n - \delta) -
\rho_s(\delta) > c(n) - 2\varepsilon$.
\end{proof}

So on one orbit, beyond the bounds of Section~\ref{sec:band}, the only
constraints are monotonicity, the triangle inequality up to the slack
of $\rho$, and condition (S) on ascents. The pairs of descent and
ascent defects are constrained only by Propositions~\ref{prop:fee}(c)
and~\ref{prop:afee}(a) and Theorem~\ref{thm:vanish}.

\begin{theorem}[Descent and ascent defects on an orbit]\label{thm:pairs}
Let $c, A : \Z \to [0, \infty)$ be nondecreasing, with $A \le c$ and
$c(n) \to 0$ as $n \to -\infty$. Then
\[
\tau(n, m) =
\begin{cases}
c(n) & \text{if } m < n,\\
0 & \text{if } m = n,\\
A(2n - m + 1) & \text{if } m > n
\end{cases}
\]
is an admissible toll over $\rho$ that satisfies the triangle
inequality, with descent defect $c$ and ascent defect $A$. So the pairs
$(c, A)$ that occur on an orbit are exactly those with $0 \le A \le c$,
both nondecreasing, and $c$ tending to $0$ in the past.
\end{theorem}

\begin{proof}
For the triangle inequality, consider $\tau(n, k)$ and a route through
$m$; routes with $m = n$ or $m = k$ are trivial. If $k < n$, the value
is $c(n)$, and the route costs at least $c(n)$: if $m < n$ its first
step is a descent from $n$, and if $m > n$ its second step is a
descent from $m$, of cost $c(m) \ge c(n)$. If $k > n$, the value is
$A(2n - k + 1) \le A(n) \le c(n)$, since $2n - k + 1 \le n$; a route
through $m < n$ or $m > k$ contains a descent from $n$ or from $m >
n$, of cost at least $c(n)$, and a route through $n < m < k$ costs at
least $A(2n - m + 1) \ge A(2n - k + 1)$. Shifting does not decrease
$c(n)$ or $A(2n - m + 1)$, which gives (b), and both tend to $0$ in the
past, which gives (c). The descent defect is $c$, and the ascent defect
is $\sup_{j \ge 1} A(n - j + 1) = A(n)$. Condition (V) holds, and so
does (S): for $n_1 < n_2 < \cdots$ and fixed $i$, $\tau(n_i, n_j) =
A(2n_i - n_j + 1) \to 0$ as $j \to \infty$. By
Corollary~\ref{cor:chains}, $\tau$ is admissible. The last statement
follows with Propositions~\ref{prop:fee}(c) and~\ref{prop:afee}(a)
and Theorem~\ref{thm:vanish}.
\end{proof}

The chain condition of Proposition~\ref{prop:afee}(b) adds nothing
here: since $c$ tends to $0$ in the past (Theorem~\ref{thm:vanish}),
a chain with $t_1$ large enough satisfies it trivially.

\begin{remark}[The spirals as extremes]\label{rem:envelopes}
Let $\tau$ be an admissible toll over $\rho$ with descent defect $c$
and no defect on ascents. Then $\tau$ lies below the source toll
$c(n)\,\mathbf{1}[m < n]$, which is admissible by Theorem~\ref{thm:pairs} with
$A = 0$; and if $\tau$ satisfies the triangle inequality, it lies above
$\tau_c$, which is admissible by Theorem~\ref{thm:occur}. For the step
$c = \kappa\,\mathbf{1}[n \ge 0]$, both are the toll
$\kappa\,\mathbf{1}[m < n \text{ and } n \ge 0]$ (for $\tau_c$ with
$n > m \ge 0$, take $k = n$), which over the arc metric is the toll
spiral. So this toll is the only admissible toll with this descent
defect, without ascent defect, that satisfies the triangle inequality.
For $c(n) = 2^n$, $\tau_c(n, m) = (2^n - 2^m)^+$ is the height spiral.
The ascent spiral charges the same $\tau_c$ on descents and $A(2n - m +
1)$ on ascents, with $A(n) = 2^{n-1} = \tau_c(n, n - 1)$; this is not
the toll of Theorem~\ref{thm:pairs}, which charges $c(n)$ on descents.
For the step and the arc metric, the pointwise lower bound
$\bigl(\tau_c(n, m) - \rho_s(n - m)\bigr)^+$ of
Section~\ref{sec:band} is $(\kappa - 2\rho)^+$ on descents from
recent sources, that is, the thin spiral. So for a step
the lower bound is itself attained by one space, and every admissible
toll with this descent defect lies, on descents, between the thin and
the toll spirals. For two steps this fails: for $c = \kappa\,\mathbf{1}[n \ge
0] + \kappa\,\mathbf{1}[n \ge 1]$ with $\kappa \ge 1$, and $\theta$ at distance
more than $1/4$ from the nearest integer, so that $\rho(0, 2) <
2\rho(0, 1)$, the triangle inequality for $\rho + \bigl(\tau_c(n, m) - \rho_s(n -
m)\bigr)^+$
along $1 \to 0 \to -1$ would read $2\kappa - \rho(0, 2) \le 2\kappa -
2\rho(0, 1)$, which is false.
\end{remark}

\section{Summary}\label{sec:summary}

Table~\ref{tab:summary} lists what is forced, and what is not, for a
bijective non-contractive map $F$ of a quasi-pseudo-metric space
satisfying (i) and a strict pair $(a, b)$.

\begin{table}[tbp]
\centering
\small
\setlength{\tabcolsep}{4pt}
\begin{tabular}{@{}lll@{}}
forced & statement & where\\
\hline
infinite strict orbit & $G^p a \ne a$; the orbit of $a$ is itself non-isometric &
Cor.~\ref{cor:orbit}\\
source carries the defect & $\D(a, G^\Delta a) > \D(a, b) - \varepsilon$
along a chain & Cor.~\ref{cor:uniform}\\
compact completed orbits & $\lams$ totally bounded on every orbit &
Thms.~\ref{thm:shadow}, \ref{thm:shadtop}\\
source bound & $\D(a, y) \le c(a)$, $c(a) \le c(Fa)$ &
Prop.~\ref{prop:fee}\\
two-sided bound & $c(G^k x) - c(G^k y) - \lams(x, y) \le \D(x, y) \le c(x)$ &
Thm.~\ref{thm:band}\\
$c \to 0$ & $c(G^k a) \to 0$; asymptotic isometry &
Thm.~\ref{thm:vanish}, Cor.~\ref{cor:rigid}\\
no host with limits & every $R_0$ host lacks upper and lower limits &
Cor.~\ref{cor:host}\\[2pt]
\hline
not forced & example & where\\
\hline
toll by source only & height and thin spirals &
Props.~\ref{prop:height}, \ref{prop:thin}\\
\begin{tabular}[t]{@{}l@{}}defect on every descent\\ from $a$ if $c(a) > 0$\end{tabular} & thin spiral ($\kappa < 1$) & Prop.~\ref{prop:thin}\\
defect on descents only & ascent spiral & Props.~\ref{prop:ascent}, \ref{prop:bilip}\\
a particular $c$ & every $c$ allowed occurs, over every $\lam$ &
Thm.~\ref{thm:occur}\\
a rate of decay & zone tolls reach $0$; $c(n) = 2^n$ does not &
Thm.~\ref{thm:zone}; Rem.~\ref{rem:rate}\\
$A$ given $c$ & every nondecreasing $A \le c$ occurs on an orbit &
Thm.~\ref{thm:pairs}\\
a particular completed orbit & every compact metrizable monothetic group &
Prop.~\ref{prop:realize}\\[2pt]
\hline
\multicolumn{3}{@{}l@{}}{the main forced properties together: $\lam$ and the two-sided
bound for $\D$ \hfill Thm.~\ref{thm:structure}}\\
\multicolumn{3}{@{}l@{}}{on one orbit: monotonicity, the triangle
inequality up to the slack, and (S) \hfill Thm.~\ref{thm:orbit}}
\end{tabular}
\caption{What is forced and what is not.}\label{tab:summary}
\end{table}
Both bounds are attained, the upper by the toll spiral and the lower
by the thin spiral, and every mixture of these two spirals has its
defect between them (Remark~\ref{rem:mix}).

The central problem of the plasticity literature, whether the unit
ball of every Banach space is plastic \cite{CKOW2016}, has a
quasi-metric form \cite[Question~7.4]{DuzerC}. Call the unit ball
$B_p = \{p \le 1\}$ of an asymmetric normed space $(X, p)$
\cite{Cobzas2013} plastic if every non-contractive bijection of
$(B_p, d_p)$, $d_p(x, y) = p(y - x)$, is an isometry. In finite
dimension, the symmetric ball $B_p \cap (-B_p)$ and sphere are rigid
for all non-contractive self-maps, for every asymmetric norm, and so
are the forward ball $B_p$ and sphere when $p$ is positive off $0$
\cite[Proposition~7.1]{DuzerC}; for
$B_p$ the positivity of $p$ is needed \cite[Example~7.2]{DuzerC}. In infinite dimension the
symmetric ball and sphere are not plastic in general
\cite[Example~7.3]{DuzerC}, and for $T_1$ asymmetric norms the
forward ball and sphere are the open cases
\cite[Question~7.4]{DuzerC}.

\section*{Data availability}

{\emergencystretch=3em
The results are machine-checked in Lean~4 using Mathlib, in
the files \texttt{Register.lean}, \texttt{Structure.lean} and
\texttt{Converses.lean}
archived with this paper at Zenodo
(\href{https://doi.org/10.5281/zenodo.24289545}{doi:\nolinkurl{10.5281/zenodo.24289545}}) together with copies of the files
\texttt{UpperLimits.lean},
\texttt{LowerLimits.lean} and \texttt{TollExamples.lean} of
\cite{DuzerC}, on which they build, with the same Lean and
Mathlib versions (Lean~4.32.0 and the Mathlib release tagged
\texttt{v4.32.0}).
The same files accompany the arXiv version
(\href{https://arxiv.org/abs/2611.03681}{arXiv:2611.03681}) as
ancillary files. In the Lean names, the asymptotic distance is
called \texttt{lam}, the descent defect \texttt{fee} and the ascent
defect \texttt{afee}, and some names keep earlier terminology (for
example \texttt{source\_pays}, \texttt{band\_lower},
\texttt{orbit\_shadow} and
\texttt{height\_\allowbreak register\_\allowbreak fee}). The files contain,
with the paper's numbering:
Theorem~\ref{thm:source} (\texttt{target\_returns},
\texttt{source\_pays}); Corollary~\ref{cor:orbit}
(\texttt{source\_\allowbreak orbit\_\allowbreak strict}, \texttt{source\_\allowbreak aperiodic});
Lemma~\ref{lem:register} (\texttt{lam\_invariant},
\texttt{lam\_tri}, \texttt{lam\_le}, \texttt{lam\_nonneg},
\texttt{defect\_mono}, \texttt{defect\_tri}, \texttt{bdd\_from});
Proposition~\ref{prop:fee} (\texttt{fee\_approx},
\texttt{defect\_le\_fee}, \texttt{fee\_mono});
Theorem~\ref{thm:band} (\texttt{band\_lower},
\texttt{band\_lower\_past}); Proposition~\ref{prop:afee}(a),(b)
(\texttt{afee\_le\_fee}, \texttt{afee\_mono},
\texttt{chain\_defect}, \texttt{fee\_le\_afee});
Proposition~\ref{prop:deep} (\texttt{deep\_\allowbreak cheap\_\allowbreak descents});
Theorem~\ref{thm:shadow} (\texttt{invariant\_\allowbreak symm\_\allowbreak precompact},
\texttt{orbit\_shadow}); the three spirals over the circle with an
irrational rotation (\texttt{height\_spiral}, \texttt{thin\_spiral},
\texttt{ascent\_spiral}), each with the quasi-pseudo-metric and $T_1$
axioms, hypothesis (i), non-contractivity of the shift, and its
strictness (on every descent for the height spiral, on every ascent
for the ascent spiral, and on the pairs of
Proposition~\ref{prop:thin}(a) for the thin spiral), and, for the
height spiral, the dependence of the toll on the
target, $\tau(1, 0) \ne \tau(1, -1)$; the asymptotic distances of the three spirals, shown to equal the
base $\rho$, and their descent defects, and the ascent defect $A(n) = 2^{n - 1}$
of the ascent spiral, computed
(\texttt{height\_\allowbreak register\_\allowbreak fee},
\texttt{thin\_\allowbreak register\_\allowbreak fee},
\texttt{ascent\_\allowbreak register\_\allowbreak fee}), with the
thin spiral attaining the lower bound of Theorem~\ref{thm:band}, for every $\kappa > 0$,
on a backward iterate with defect strictly between $0$ and the descent defect
(\texttt{thin\_band\_attained}); the bundled \texttt{thin\_spiral}
and \texttt{thin\_\allowbreak register\_\allowbreak fee} appear in the
axiom lists at the instance $\kappa = 2\rho(0, -1)$, while the general-$\kappa$ statements
are the compiled lemmas \texttt{thin\_isQPM}, \texttt{thin\_t1},
\texttt{thin\_seqHP}, \texttt{thin\_shift\_nc}, \texttt{close\_return},
\texttt{thin\_paid} and \texttt{thin\_shift\_strict}
(Proposition~\ref{prop:thin}(a)), \texttt{thin\_free}
(Proposition~\ref{prop:thin}(b)),
\texttt{thin\_lam}, \texttt{thin\_defect} (the defect
$(\kappa - 2\rho)^+$ of Proposition~\ref{prop:thin} on descents from
recent sources), \texttt{thin\_fee\_recent}
and \texttt{thin\_fee\_deep};
Corollary~\ref{cor:uniform} (\texttt{uniform\_monochromy});
the equality of the backward limits of
Proposition~\ref{prop:afee}(c) in its infimum form
(\texttt{fee\_afee\_inf}; by \texttt{fee\_iter} and
\texttt{afee\_iter} the infima are the limits); and
Theorem~\ref{thm:zone}
(\texttt{addSubmonoid\_\allowbreak nat\_\allowbreak chain\_\allowbreak stabilizes},
\texttt{zone\_stabilizes}, \texttt{zone\_fee\_zero}; the constant
ladder is the compiled lemma
\texttt{zd\_\allowbreak ladder\_\allowbreak const}), with the
asymptotic distance and descent defect as defined there; and, in
\texttt{Structure.lean}, \texttt{lam\_self}, which completes
Lemma~\ref{lem:register}(a), Lemma~\ref{lem:feetarget}
(\texttt{fee\_target}), Theorem~\ref{thm:vanish}
(\texttt{fee\_vanishes}), Corollary~\ref{cor:rigid}
(\texttt{deep\_\allowbreak sources\_\allowbreak rigid},
\texttt{afee\_vanishes}) and Theorem~\ref{thm:structure}(a),(b)
(\texttt{structure\_theorem}).
\texttt{Converses.lean}, new in version 3, contains: the maximality
of the asymptotic distance, the last sentence of
Remark~\ref{rem:witness} (\texttt{lam\_greatest},
\texttt{le\_\allowbreak lam\_\allowbreak of\_\allowbreak iso});
Proposition~\ref{prop:bilip}, stated through the pulled-back
distance and using only the upper Lipschitz bound, with the
conjugacy not assumed bijective or equivariant
(\texttt{ascent\_\allowbreak not\_\allowbreak tolled},
\texttt{ascent\_\allowbreak spiral\_\allowbreak not\_\allowbreak tolled});
the statements of Proposition~\ref{prop:realize} about $(\Z, d)$
(\texttt{every\_\allowbreak completion\_\allowbreak occurs}; for
finite $K$, \texttt{dR\_isQPM} and \texttt{dR\_seqHP}),
with the $T_1$ argument of Theorem~\ref{thm:shadtop}
(\texttt{qpm\_t1}); and Theorem~\ref{thm:pairs}, in both
directions (\texttt{pairs\_occur}, \texttt{pairs\_necessary});
condition (i) for the new distances is proved directly, without the
Ramsey theorem of Section~\ref{sec:occur}. Each of the three files
ends with a list of \texttt{\#print axioms} commands, $39$ in
all. Every
theorem in these lists depends only on the axioms \texttt{propext},
\texttt{Classical.choice} and \texttt{Quot.sound}. The supporting
lemmas named above are proved in the same files but do not appear in
these lists.

The following are not mechanized:
\begin{itemize}
\item the asymptotic distance and descent defect of the toll spiral of \cite{Duzer2026},
and the fact that it attains the upper bound of
Theorem~\ref{thm:band};
\item Theorem~\ref{thm:shadtop}, the completion of each orbit to
a compact group in Section~\ref{sec:shadow}, and, in
Proposition~\ref{prop:realize}, the identification of the completed
orbit with $(K, \ell)$, with the paragraph after it;
\item the statements about the sets $U_\varepsilon$ at the end of
Theorem~\ref{thm:structure}, each one line from the conjuncts of
\texttt{structure\_theorem}, and the isometry form in
Corollary~\ref{cor:rigid}, one line from
\texttt{deep\_\allowbreak sources\_\allowbreak rigid},
\texttt{lam\_invariant} and \texttt{lam\_le};
\item the last two sentences of Corollary~\ref{cor:orbit}(a), and
the normal form described after Corollary~\ref{cor:uniform};
\item the sentences beginning ``In particular'' in
Theorem~\ref{thm:source} and Corollary~\ref{cor:uniform}, whose
mechanized forms (\texttt{source\_pays},
\texttt{uniform\_\allowbreak monochromy}) take the hypothesis on
$\eta$ and $M$ directly;
\item the statement in Section~\ref{sec:depth} that $\D = 0$ exactly
when $F$ is an isometry;
\item the passage from the sequential form of total boundedness, which
is what \texttt{invariant\_\allowbreak symm\_\allowbreak precompact}
and \texttt{orbit\_shadow} state, to total boundedness itself in
Theorem~\ref{thm:structure}(a) and in the paragraph after
Theorem~\ref{thm:shadow};
\item Corollary~\ref{cor:tollanswer}, which combines the mechanized
parts of Theorem~\ref{thm:structure} with these results on the
completed orbits;
\item the statement after Proposition~\ref{prop:thin} that targets as
in (b) exist exactly when $\kappa < 1$, and the comparison with
Theorem~\ref{thm:band} after Proposition~\ref{prop:extremal};
\item Remarks~\ref{rem:asym}, \ref{rem:readings}, \ref{rem:mix}
and~\ref{rem:z1z2}, Remark~\ref{rem:witness} except its last
sentence, and Propositions~\ref{prop:extremal}
and~\ref{prop:tolled} (with the two paragraphs after it);
\item the discussion after Remark~\ref{rem:readings}: the heuristic
of \cite[Remark~5.8]{Duzer2026}, the integer layer of
\cite[Proposition~6.1]{DuzerC}, and the unboundedness of every
nonzero toll of the form $(h(n) - h(m))^+$;
\item the total boundedness of backward orbits for
$d(x, y) + d(y, x)$ and the near-returns on their tails, noted after
Corollary~\ref{cor:rigid}, and the toll-spiral example that follows
them;
\item the reformulation of Theorem~\ref{thm:band} through
$\tau_c$ after it, Corollary~\ref{cor:host} and the paragraph after
it, and Section~\ref{sec:occur} except Theorem~\ref{thm:pairs};
\item Remark~\ref{rem:rate}, except for the descent defects $c(n) = 2^n$ of the
height and ascent spirals, which are the mechanized results
\texttt{height\_\allowbreak register\_\allowbreak fee} and
\texttt{ascent\_\allowbreak register\_\allowbreak fee}; the
conclusion that Theorem~\ref{thm:vanish} cannot be improved is argued
on paper.
\end{itemize}
The files of \cite{DuzerC} are included, and compiled, because
the present files import lemmas from them. Results otherwise
cited from other papers, such as the equivalence of (i) with
hereditary precompactness \cite[Theorem~3]{KMRV1993} and the facts
on unit balls in Section~\ref{sec:summary}, are not mechanized.

\par\medskip\noindent
\begin{minipage}{\linewidth}
\hspace*{17pt}The manifest \texttt{ArtifactHashes.txt} records the
SHA-256 hash of each Lean file and of the three toolchain pins, and its own
SHA-256 hash is
\begin{center}\small
\texttt{b2b7d7816c241d73946aef6da5ac9433a95d7d7c74496178698816bc3e48a0ed}
\end{center}
\end{minipage}\par}

\section*{Acknowledgements}

The author developed the mathematics in a draft in French and
used AI to translate it into a
first English version. The author
corrected this version wherever it departed from the French and
rewrote it, keeping some of its standard phrasing. The author
then used AI to check the grammar and reviewed the changes.

\end{document}